\documentclass[reqno,onefignum,onetabnum]{siamart171218}
\usepackage{subdepth}
\usepackage{overpic}
\usepackage{amsmath,amstext,amsbsy,amssymb,mathdots}
\usepackage{booktabs,bm}
\usepackage{mathrsfs}
\usepackage{tikz,cite}
\usetikzlibrary{intersections}
\usepackage[export]{adjustbox}

\usepackage{courier}
\usetikzlibrary{positioning}
\usetikzlibrary{shapes,arrows}
\usepackage[fleqn,tbtags]{mathtools}
\usepackage{amsfonts}
\usepackage{relsize}
\usepackage{xcolor,colortbl}
\usepackage{psfrag}
\usepackage{alltt}
\usepackage{arydshln}
\usepackage{enumerate}
\usepackage{epstopdf}
\usepackage{enumitem}
\usepackage{multirow}
\usepackage{algorithmic}
\usepackage{framed}
\usepackage{cancel}
\usepackage{tikz}
\usepackage{tkz-euclide}
\usepackage{graphicx} 
\graphicspath{{images/}}
\usepackage[caption=false]{subfig}
\usepackage{pgfplots}

\usetikzlibrary{patterns}

\usepackage{appendix}

\newcommand*{\CopyCounter}[2]{%
  \expandafter\def\csname c@#2\endcsname{\csname c@#1\endcsname}%
  \expandafter\def\csname p@#2\endcsname{\csname p@#1\endcsname}%
  \expandafter\def\csname the#2\endcsname{\csname the#1\endcsname}}

\CopyCounter{Theorem}{ProposedProblem}
\CopyCounter{Theorem}{Proposition}
\CopyCounter{Theorem}{Property}
\CopyCounter{Theorem}{Claim}
\CopyCounter{Theorem}{Lemma}
\CopyCounter{Theorem}{Corollary}
\CopyCounter{Theorem}{Conjecture}
\CopyCounter{Theorem}{Definition}
\CopyCounter{Theorem}{Example}
\CopyCounter{Theorem}{Remark}
\CopyCounter{Theorem}{Question}
\CopyCounter{Theorem}{Condition}
\CopyCounter{Theorem}{Criterion}
\CopyCounter{Theorem}{Observation}
\theoremstyle{plain}

\newcommand{\domain}{\Omega}

\newcommand{\boundary}{\partial \domain}

\newcommand{\xtilde}{\bm{\tilde{x}}}

\newcommand{\ptilde}{\tilde{p}}

\newcommand{\xhat}{\bm{\hat{x}}}

\newcommand{\Xhat}{\hat{X}}

\newcommand{\x}{\bm{x}}

\newcommand{\ba}{\bm{a}}

\newcommand{\Bf}{\bm{f}}

\newcommand{\bs}{\bm{s}}

\newcommand{\y}{\bm{y}}
\newcommand{\z}{\bm{z}}

\newcommand{\J}{{\cal{J}}}

\newcommand{\R}{\mathbb{R}}
\newcommand{\E}{\mathbb{E}}

\newcommand{\G}{\mathcal{G}}
\newcommand{\M}{\mathcal{M}}

\newcommand{\V}{\mathcal{V}}

\newcommand{\qtilde}{\tilde{q}}

\newcommand{\qbar}{\bar{q}}
\newcommand{\qbeta}{q_{\beta}^{}}

\newcommand{\source}{\bm{\x_0^{}}}
\newcommand{\sa}{\bs}
\newcommand{\sw}{\bm{\bar{s}}}
\newcommand{\srs}{\bm{\bs_{\beta}^{}}}
\newcommand{\saw}{\bm{\bs_{c}^{}}}

\newcommand{\PP}{\mathbb{P}}

\newcommand{\domainT}{\domain^{}_T}
\newcommand{\omegat}{\domainT}
\newcommand{\worst}{C}
\newcommand{\omegatc}{\Omega_{T,C}^{}}
\newcommand{\omegac}{\Omega_{C}^{}}
\newcommand{\phat}{\hat{p}}
\newcommand{\measure}{\hat{\mu}}

\newcommand{\bptilde}{\bm{\tilde{p}}}
\newcommand{\bstilde}{\bm{\tilde{s}}}

\newcommand{\bphat}{\bm{\hat{p}}}

\DeclareMathOperator*{\argmin}{arg\,min}

\usepackage{xargs}
\usepackage{comment}
\usepackage[colorinlistoftodos,prependcaption,textsize=tiny]{todonotes}
\reversemarginpar
\newcommandx{\avtodo}[2][1=]{\todo[linecolor=red,backgroundcolor=red!25,bordercolor=red,#1]{#2}}
\newcommandx{\mgtodo}[2][1=]{\todo[linecolor=blue,backgroundcolor=blue!25,bordercolor=blue,#1]{#2}}

\colorlet{lightgray}{gray!40}

\newcolumntype{"}{@{\hskip\tabcolsep\vrule width 1pt\hskip\tabcolsep}}
\setlist[description]{font=\normalfont\space}

\usepackage{placeins}

\usepackage{lipsum}
\ifpdf
  \DeclareGraphicsExtensions{.eps,.pdf,.png,.jpg}
\else
  \DeclareGraphicsExtensions{.eps}
\fi

\newsiamremark{remark}{Remark}
\newsiamremark{hypothesis}{Hypothesis}
\crefname{hypothesis}{Hypothesis}{Hypotheses}

\newif\ifshowchanges
\showchangesfalse

\ifshowchanges
\newcommand{\change}[1]{{\color{blue} #1}}
\else
\newcommand{\change}[1]{{\color{.} #1}}
\fi

\headers{Path-Planning Under Initial Uncertainty}{Qi, Dhillon, and Vladimirsky}

\title{Optimality and robustness in path-planning under initial uncertainty
}
\author{Dongping Qi\thanks{Center for Applied Mathematics, Cornell University, Ithaca, NY 14853
  (\email{dq48@cornell.edu}).}
\and Adam Dhillon\thanks{Department of Mathematics, University of California, Berkeley, CA 94720
  (\email{adamdhillon@berkeley.edu}).}   
\and Alexander Vladimirsky\thanks{Department of Mathematics and Center for Applied Mathematics, Cornell University, Ithaca, NY 14853
  (\email{vladimirsky@cornell.edu}).}
}

\usepackage{amsopn}

\begin{document}
\maketitle


\begin{abstract}
Classical deterministic optimal control problems assume full information about the controlled process.  The theory of control for general partially-observable processes is  powerful, but the methods are computationally expensive and typically address the problems with stochastic dynamics and continuous (directly unobserved) stochastic perturbations.  In this paper we focus on path planning problems which are in between -- deterministic, but with an initial uncertainty on either the target or the running cost on parts of the domain.  That uncertainty is later removed at some time $T$, and the goal is to choose the optimal trajectory until then. 
We address this challenge for three different models of information acquisition: with fixed $T$, discretely distributed and exponentially distributed random $T$. 
We develop models and numerical methods suitable for multiple notions of optimality: based on the average-case performance, the worst-case performance, the average constrained by the worst, the average performance with probabilistic constraints on the bad outcomes, risk-sensitivity, and distributional-robustness.
We illustrate our approach using examples of pursuing random targets identified at a (possibly random) later time $T$.
\end{abstract}




\section{Introduction}
\label{sec:Intro}

A common task in robotic navigation is to find a trajectory that minimizes cumulative cost incurred on the way to a target.
If the system state is perfectly known in real time, this can be accomplished in the framework of dynamic programming by solving
Hamilton-Jacobi-Bellman (HJB) partial differential equations (PDEs) \cite{Bardi_book_original, FlemingSoner}.
The task of controlling {\em partially-observable} processes is generally harder due to \change{the} added challenge of estimating the process state probabilistically
\cite{zarchan2013fundamentals, kaelbling1998planning}.  
Whether one aims for on-average optimality or robustness, the methods in this general space are focused on the uncertainty due to 
frequent stochastic perturbations in between observations.
In contrast, we are interested in controlling processes with uncertainty due to structured or delayed information acquisition patterns. 
Here we focus on the simplest subset of such problem, in which everything is  
deterministic, but some important aspects
of the process are a priori unknown.   
We start with a probabilistic description of the global environment, dynamics, costs or targets,
and then need to control the system until this initial uncertainty is resolved.

In \change{a} discrete setting, a classical example of 
such monotonically non-increasing
uncertainty is the ``Canadian Traveller Problem'' \cite{PAPADIMITRIOU1991127}, in which one aims to find a cost-minimizing
path to a target on a partially known graph: some of the graph edges might be absent
with pre-specified/known probabilities, but whether they actually are present is only learned later, as we enter nodes adjacent to them.  
A generalization allowing random  
edge costs and more complicated information-gathering on a graph is known as a  ``Stochastic Shortest Path with Recourse''  
\cite{polychronopoulos1996stochastic}.  While dynamic programming algorithms are well known for these problems, they operate on a much larger information-enriched state space, 
and their computational cost is typically prohibitive.

In this paper we examine path-planning problems in a continuous setting but with a fairly simple information acquisition pattern:
the full information about the system and the environment is instantaneously gained at some ``certainty time'' $T$, and the problem remains deterministic  
thereafter.
This simplification allows us to develop efficient algorithms for optimizing the pre-certainty control and also helps in comparing the merits of several notions of robustness.
For motivation, consider a semi-autonomous rover exploring a sequence of locations on Mars with infrequent guidance from a command center on Earth.
The center might identify a set of possible future targets (with \change{a} probability specified for each of them a priori),  
but the final choice will be identified at a later time $T$.   While the rover could wait in place until then, its cumulative cost-to-target can be reduced by 
choosing a suitable waypoint and starting to move to it already now.
To simplify the exposition,
we describe all algorithms in the context of such  {\em initial target uncertainty} problems in 2D and with isotropic cost/dynamics.  
But our general approach is broader and also covers other types of initial uncertainty. 
In Appendix \ref{sec:storm_fronts}, we show how it can be applied to a simple environment-uncertainty problem: planning a path for an airplane to avoid a storm front whose position remains uncertain until \change{the} time $T$.

In section \ref{sec:Determ},
we start by reviewing the standard deterministic setting: planning an optimal path to a known target by solving a single HJB equation.
We then introduce the target uncertainty in section \ref{sec:AverageCase} and develop methods for optimizing the {\em average-case performance}.
Interestingly, the methods for doing this depend entirely on what is known about the certainty time $T$.  Once the min-cost to each potential target is known,
we show that, for a fixed/constant $T,$ the expected-cost-minimizing waypoint is found by optimizing over a constraint set determined by solving a single stationary HJB equation. 
{
}
If $T$ is random and has a known discrete distribution, we show that optimal waypoints are found by solving a sequence of time-dependent
HJB equations.  
For simplicity, our method is described with the likelihood of individual targets viewed as independent of the realized value of $T$,
but we later show in Appendix \ref{sec:drones_example} that this assumption is not essential.
If $T$ is an exponentially distributed random variable, finding an optimal waypoint is not enough -- 
we also need to find 
an optimal path toward that waypoint
since the target might be identified before we reach it.  
We explain how this can be accomplished by solving a quasi-variational HJB inequality of obstacle type.

The issue of robustness under initial uncertainty is the subject of section \ref{sec:Robustness}.
Robust control is an important and well-developed research area, with $H^\infty$ methods \cite{bacsar2008h} particularly popular in stochastic systems subject to frequent or continuous random disturbances.  
But in our setting, this approach is not directly applicable, though one can still mitigate the risk of ``unlucky scenarios'' and guard against modeling errors.   
Focusing primarily on the ``fixed $T$'' case, we develop methods for several competing notions of robustness:
optimizing the worst\change{-}case scenario,  risk-sensitive optimization, and optimizing the average case performance subject to constraints (either hard or probabilistic) on 
bad outcomes.  
The case of probabilistic constraints is particularly interesting since the optimal policy for selecting waypoints turns out to be non-deterministic.  
We find a geometric interpretation of this challenge and leverage it to construct an efficient algorithm in subsection \ref{ss:probabilistic}. 
We further examine the Distributionally Robust Optimization (DRO) approach to guard against mistakes in \change{the} perceived relative likelihood of potential targets. 
In addition, we also provide a sub-optimality bound for the case when the set of targets is either misidentified or intentionally subsampled.

Many of the problems described above are reminiscent of those handled in discrete setting\change{s} by the techniques of {\em multistage stochastic programming} \cite{shapiro2021lectures}.  However, their naive use would be prohibitive here since our state and action spaces are infinite/continuous.
In addition, our approach based on dynamic programming allows solving these multi-stage optimization problems for all initial configurations simultaneously.
These distinctions (and the computational efficiency gained by exploiting the problem structure) are highlighted throughout the paper.
We conclude by discussing possible generalizations of our approach and directions for future work in section \ref{sec:Conclusions}.


\section{Deterministic setting overview}
\label{sec:Determ}

We first consider a fully deterministic controlled process with a time-dependent process state $\y \in 
\R^d$
and a terminal time $T.$
This is the setting of the classical optimal control theory and we include only a brief overview,
referring to  \cite{Bardi_book_original} for technical details.

We will assume that the process dynamics is defined by an ODE
\begin{equation}
\label{eq:general_dynamics}
\begin{cases}
\y'(s)  = \Bf \left(\y(s), \ba(s), s \right),  \qquad s \in [t, T];\\
\y(t)  = \x \in  
\R^d. 
\end{cases}
\end{equation}
The controller has a compact set $A$ of available control values and chooses a measurable control function
$\ba : \R \mapsto A$ to guide the evolution of the process state. 
The ``velocity function'' $\Bf:(\R^d \times A \times \R) \mapsto \R^d$, the initial state and time $(\x,t),$ and the terminal time $T$ are all explicitly known ahead of time.  Given a running cost $K:\R^d \times A \times \R \mapsto \R$ and a terminal cost $q : \R^d \mapsto \R,$ one can define the cost corresponding to a control $\ba(\cdot)$ starting from $(\x,t)$ as
\begin{equation}
\label{eq:objective}
\J\left( \x,t, \ba(\cdot) \right) \, = \, \int_t^T K \left( \y(s), \ba(s), s \right) \, ds \, + \, q(\y(T)).
\end{equation}
A {\em value  function} is then defined as the minimal cost-to-termination from a given location; i.e.,
$
u(\x,t) \; = \; \inf_{\ba(\cdot)} \J \left( \x,t,  \ba(\cdot) \right).
$
The principle of dynamic programming is based on a {\em tail optimality} property of optimal controls; i.e., 
\[
u(\x,t) \, = \, 
\min\limits_{\ba(\cdot)}
\left\{
\int_t^{t+\tau} K \left( \y(s), \ba(s), s \right) \, ds \, + \, u\left(\y(t+\tau), t+\tau \right)
\right\}
\]
should hold for any small $\tau>0$.   Assuming that the value function is smooth, one can use a Taylor series expansion and then let $\tau \to 0$ to derive a time-dependent Hamilton-Jacobi PDE that $u$ must satisfy:
\begin{equation}
\label{eq:tHJB}
u_t (\x, t)
\, + \,
\min\limits_{\ba \in A}
\left\{
K(\x, \ba, t) + \nabla u (\x, t) \cdot \Bf ( \x, \ba, t )
\right\}
\, = \, 0,
\quad \forall (\x,t) \in \R^d \times [0,T)
\end{equation}
subject to a {\it terminal condition} $u(\x,T)=q(\x), \; \forall \x \in \R^d$.
Since the value function is generally not smooth, this PDE usually  does not have a classical solution.
The non-uniqueness of weak (locally Lipschitz) solutions made it necessary to introduce additional test conditions
\cite{crandall1983viscosity} that pick out a {\em viscosity solution} of this PDE -- the unique weak solution that coincides with the value function of the original control problem.

Even though the approach presented in this paper for treating the initial uncertainty is general, all of our examples will be based on {\it isotropic} problems,  in which the cost and speed of motion depend only on $\x$.
We let $A=\{ \ba \in \R^d \, \mid \, |\ba| = 1 \}$ and interpret $\ba$ as the chosen direction of motion.  Then
$K(\x, \ba, t ) = K(\x,t)$ and $\Bf ( \x, \ba, t ) = f (\x,t) \ba$, with $f$ encoding the speed of motion through the point $\x$.
In this case, the optimal direction is known analytically: $\ba^* = -\nabla u / |\nabla u|$ and \eqref{eq:tHJB} 
reduces to a {\it time-dependent Eikonal equation}
\begin{equation}
\label{eq:tEik}
u_t \, - \, | \nabla u | f(\x,t) \, + \,  K(\x,t) \; = \; 0;
\qquad \forall \x \in \R^d, \, t \in [0,T).
\end{equation}
Both the isotropic and the general problems can be  
similarly posed on a compact domain $\domain \subset \R^d$,
with additional boundary conditions $u=g$ posed on $\boundary \times [0,T]$, which can be interpreted as a cost 
of stopping the process prematurely as soon as it exits $\domain$ at some $t < T$.  
Setting $g=+\infty$ can be used to restrict the process to this domain\footnote{
To define such  {\em domain-constrained viscosity solutions}, one has to treat the boundary conditions ``in the viscosity sense'' \cite{Bardi_book_original}
and we adopt this interpretation throughout the paper.
Intuitively, this means that the infinite exit cost $g$ is charged not for touching $\boundary$ but for attempting to leave $\domain.$
This allows traveling along $\boundary$ (including the obstacle boundaries) and, combined with our isotropy assumption, ensures the existence of optimal controls for every starting position.}
and force the trajectories to avoid any obstacles (by simply excluding them from $\domain$).
Setting $g(\x,t)=0$ on any closed {\em target set} $\Gamma \subset \boundary$ makes $u(\x,t)$ encode the minimal cost to target $\Gamma$, with $q$ specifying a penalty 
in case we cannot (or choose not to) reach
$\Gamma$ by the time $T$.  The latter case is further simplified when $K$ and $f$ do not depend on $t$ and $T=+\infty$ (i.e., there is no time-restriction for reaching the target).  The corresponding value function satisfies a (stationary) {\it Eikonal equation}
\begin{equation}
\label{eq:Eik}
| \nabla u | f(\x) \; = \;  K(\x), \qquad \forall \x \in \domain; \qquad
 u(\Gamma) = 0; \quad u(\boundary \backslash \Gamma) = +\infty.
\end{equation}
If $K \equiv 1$, $u(\x)$ simply encodes the minimum time to the target $\Gamma$ (while staying inside $\domain$)
and the gradient descent  
yields time-optimal trajectories.
When $\Gamma$ consists of only a single point $\source$, we denote the value function of this min-cost-from-$\source$-to-a-target-at-$\x$ problem as $u(\x;\source)$.

\begin{figure}[h]
\centering
$
\begin{array}{cc}
\includegraphics[width = 0.46\textwidth]{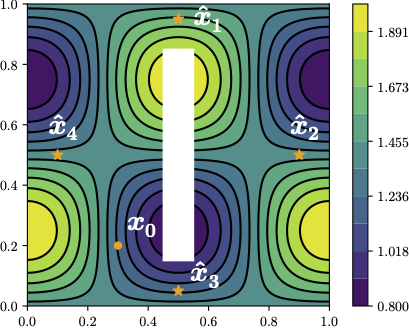} & \quad 
\includegraphics[width = 0.45\textwidth]{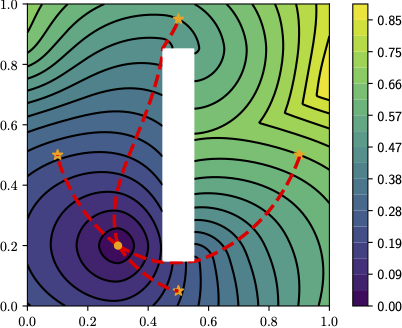}\\
(A) & (B)
\end{array}
$ 
\caption{Time-optimal path planning.
(A): a contour plot of the speed $f = 1.4+0.6\cos(2\pi x)\sin(2\pi y)$ on a unit square domain with one rectangular obstacle (shown in white) at $(0.45,0.55)\times(0.15,0.85).$
The starting point $\source = (0.3,0.2)$ is shown by an orange dot.
The approach finds optimal trajectories to all target locations, but we highlight four potential targets (shown by orange stars and numbered clockwise, starting from the top): $\xhat_1 = (0.5, 0.95), \xhat_2 = (0.9,0.5), \xhat_3 = (0.5, 0.05),$ and $\xhat_4 = (0.1, 0.5).$
(B): a contour plot of $u(\x;\source)$ and the optimal trajectories (shown by dashed red lines) to these four targets. 
}
\label{fig:determ}
\end{figure}

Another convenient property of the Eikonal equation is that an optimal trajectory from a starting position $\source$ to a target is the same (up to a time-reversal) as an optimal trajectory from that target to $\source$.  
\change{The optimal} trajectory from $\source$ to any target $\x$ can be  found by following $(-\nabla u(\x;\source))$ from $\x$ to $\source$ and then tracing it backward.  In Figure \ref{fig:determ} we present an example of this on a 2D domain with a single rectangular obstacle.  Since $K=1$, we are finding the time-optimal trajectories from $\source=(0.3,0.2)$ under the speed function $f = 1.4+0.6\cos(2\pi x)\sin(2\pi y)$ shown in Figure \ref{fig:determ}A.  In Figure \ref{fig:determ}B we show the optimal trajectories from $\source$ to 4 different targets as well as the level sets of the function $u(\x;\source)$.  In subsequent sections, this general setup will be used as our representative example to explore the initial target-uncertainty.

The viscosity solutions to the Hamilton-Jacobi-Bellman (HJB) equations are usually unavailable in analytic form and numerical approximations are thus unavoidable.
Numerical methods for time-dependent and stationary HJB equations have been an active area of research in the last 25 years, with particularly many efficient techniques developed for the stationary Eikonal equation \eqref{eq:Eik} and its anisotropic generalizations.   
The first challenge is to choose a discretization converging to the viscosity solution.  The simplest grid discretization approach (employed here) is to use the first-order {\em upwind divided differences} to approximate the partial derivatives of $u$; e.g., with $d=2$ and a space-time grid based on $(x_i, y_j, t_n) = (ih, jh, n \Delta t),$ we can 
\change{
define the standard one-sided divided differences as
$$
D^{\pm x} U^n_{i,j} \, = \, \frac{U^n_{i\pm1,j} - U^n_{i,j}}{ \pm h}
\qquad \text{and} \qquad
D^{\pm y} U^n_{i,j} \, = \, \frac{U^n_{i,j\pm1} - U^n_{i,j}}{ \pm h},
$$
then plug
}
\begin{align}
\nonumber
u_t (x_i, y_j, t_n) &\approx D^{-t} U^n_{i,j} \, = \, \frac{U^n_{i,j} - U^{n-1}_{i,j}}{\Delta t}, \\ 
\nonumber
\left[ u_x (x_i, y_j, t_n) \right]^2 &\approx \left[ \max \left\{D^{-x}U^n_{i,j}, \, -D^{+x}U^n_{i,j}, \, 0 \right\}  \right]^2 , \\
\label{eq:upwind_discr}
\left[ u_y (x_i, y_j, t_n) \right]^2 &\approx \left[ \max \left\{D^{-y}U^n_{i,j}, \, -D^{+y}U^n_{i,j}, \, 0 \right\}  \right]^2 
\end{align}
into equation \eqref{eq:tEik} and  
solve for  $U^{n-1}_{i,j}.$ (The state constraints are handled by setting $U$ values to $+\infty$ outside of $\domain$, including inside the obstacles.)

This discretization guarantees a monotone dependence on the neighboring grid values, and this monotonicity yields convergence to the viscosity solution \cite{barles1991convergence}. For the time\change{-}dependent PDE \eqref{eq:tEik}, a time-explicit discretization makes the issues of efficiency trivial: the solution is computed by time-marching, computing one time-slice at the time, from the terminal time $t=T$ to the initial time $t=0,$ with the  $O(NM)$computational cost\footnote{Throughout the paper, we will assume that PDEs are discretized on a spatial grid with a total of $M$ gridpoints and, whenever these PDEs are time-dependent, we will assume that there is a total of $N$ time-slices.  Since our time discretizations are explicit, these $N$ and $M$ are not really independent: e.g., with $d=2$, Courant-Friedrichs-Lewy stability condition yields 
$N = O(M^{1/2}).$  But we will not directly use this expression in computational complexity estimates since our goal is to emphasize the general space-time structure and other
(time-implicit \cite{VladZheng} or semi-Lagrangian \cite{falcone2013semi}) discretizations could be used to avoid this dependence.}.
But in the stationary case \eqref{eq:Eik}, the resulting system of discretized equations for $U_{i,j}$\change{s}  
is non-linear, non-smooth, and coupled, adding the challenge of solving that discretized system efficiently.  A number of methods are based on competing ideas for handling this task: from Fast Marching 
(e.g., \cite{Tsitsiklis}, \cite{sethian1996fast}), 
to Fast Sweeping (e.g., \cite{boue1999}, \cite{zhao2005fast}, \cite{tsai2003fast}),
to their hybrid versions \cite{chacon2012fast, chacon2015parallel} and  methods mirroring the logic of label-correcting algorithms on graphs (e.g., \cite{jeong2007interactive, bak2010some}).
When optimal paths are needed from a small number of starting positions only, additional gains in efficiency can be also obtained through a dynamic restriction of 
\change{the} computational domain \cite{yershov2013simplicial, ClawsonChaconVlad}.  
Many of these methods have been generalized to simplicial meshes, problems on manifolds and higher order accurate schemes.  
Other approaches have also been developed using semi-Lagrangian \cite{falcone1994minimum, falcone2013semi, FalconeFMM, potter2019ordered} and discontinuous Galerkin discretizations \cite{li2008second, zhang2011uniformly}.

For the purposes of the current paper, we are indifferent to the choice of \change{a} specific numerical method.  For the sake of simplicity, our actual implementation uses the first-order upwind discretization 
based on \eqref{eq:upwind_discr}
and the Fast Marching solver (in stationary cases),
which can find a discretized solution of \eqref{eq:Eik} in $O(M \log M)$ operations on a spatial grid with $M$ gridpoints.
But other numerical techniques could be used as well whenever we need to solve a fully deterministic problem.  
Our real focus is on treating initial uncertainties -- either in the target $\Gamma$ or in the running cost $K$ on parts of $\domain$.


\section{Average case optimality under initial uncertainty}
\label{sec:AverageCase}

We begin by introducing our main sample problem of ``optimally pursuing random targets''.
We will assume that the dynamics and the running cost are isotropic and time-independent (i.e., $K= K(\x)$ and $\Bf(\x,\ba) = f(\x) \ba$, with $|\ba|=1$), and the goal is to minimize the total cost of reaching a target $\xhat \in \domain \subset \R^2$.
When $\xhat$ is known and there is no deadline,  one can do this by solving 
\eqref{eq:Eik} with $\Gamma = \left\{ \xhat \right \}.$

But what if we only know a set of possible targets 
$\Xhat =   \left\{ \xhat_1, \ldots,  \xhat_m \right \}$ and the corresponding probability distribution
$\bphat = (\phat_1, \ldots, \phat_m)$?
We can compute each ``min cost to the i-th target'' $u(\x;\xhat_i)$ by solving  \eqref{eq:Eik} with $\Gamma = \left\{ \xhat_i \right \},$ and from now on we will use $u_i(\x)$ to denote $u(\x;\xhat_i)$ as long as there is no ambiguity. 
Then we can define the ``expected min cost to target''
\begin{equation}
\label{eq:q}
q(\x) \; = \; \E^{}_{\xhat} u(\x) = \sum\limits_{i=1}^m \phat_i u_i(\x).
\end{equation} 
Here $u(\x)$ is a random variable satisfying $P(u(\x) = u_i(\x)) = \phat_i$.
The global minima of $q$ suggest the best starting point(s) if the target is immediately revealed.   
But what if we start elsewhere and what if the actual $\xhat$ is only learned at a later time $T$?
The real question is what to do until then 
since  the gradient descent  in ``correct'' $u_i$ will define the trajectory after  $\xhat = \x_i$ is known. 


\subsection{Fixed certainty time $T$}
\label{sec:FixedT}

For a fixed and known certainty time $T$, the answer is encoded by the solution of \eqref{eq:tEik} with the terminal condition
$u = q$ on $\domain \times \{T\}$ and the boundary condition $u = +\infty$ on $\boundary \times [0,T).$ 

The case $K=1$ is particularly simple since the total cost incurred until the target identification is $T$ regardless of the  chosen control/trajectory.  Thus, the goal then becomes to minimize the ``terminal cost'' $q$ among the points reachable from our starting position $\source$.  
Solving  \eqref{eq:tEik} has the advantage of encoding the answer for all starting positions simultaneously, but for any specific $\source$ we can further simplify the computation.
In this case, we solve \eqref{eq:Eik} with $\Gamma = \{\source\}$ 
to find the min-time-from-$\source$ (denoted $u(\x;\source)$; 
see the example in Figure \ref{fig:determ}B)
and define the reachable set\footnote{
With $K=1$ and fixed $T$, if we reach $\x$ too early, we can simply wait there without affecting the cumulative cost.}
\[
\domainT(\source) = \{ \x \in \domain \, \mid \, u(\x;\source) \leq T \}.
\]
To simplify the notation, we will use $\domainT$ whenever $\source$ is clear from the context.
The optimal expected time-to-target is then obtained by traveling from $\source$ to any waypoint
\[
\sa \in \argmin\limits_{
\x \in \domainT(\source)}  q(\x).
\]

Since there are a total of $m$ targets (each with its own value function $u_i$), 
$q$ can be found on the entire $\domain$ in $O(mM\log M)$ operations.
After that, for any specific starting position $\source$ the optimal waypoint $\sa$ can be found in
$O(\tilde{M} \log \tilde{M})$ operations needed to determine the $\tilde{M} \leq M$ gridpoints falling in the set  $\domainT(\source).$

\begin{figure}[h]
\centering
$
\begin{array}{cc}
\includegraphics[width = 0.45\textwidth]{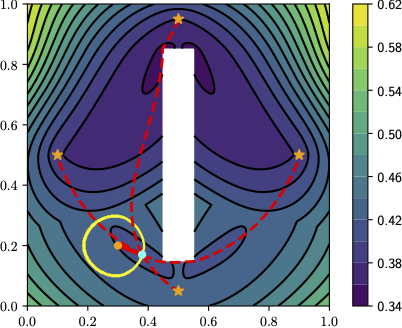} & \quad 
\includegraphics[width = 0.45\textwidth]{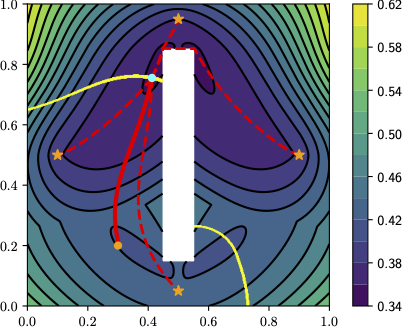}\\
(A) & (B)
\end{array}
$
\caption{Fixed certainty time.  (A): $T=0.08.$  (B): $T=0.4.$
The speed $f$, domain geometry, starting position $\source$ (orange dot)
and four targets $\xhat_1, \ldots, \xhat_4$ (orange stars) are the same as in Figure \ref{fig:determ}.
The corresponding target probabilities (clockwise, starting at the top) are $\bm{\phat} = (0.2, 0.3, 0.2, 0.3).$
Boundaries of $\domainT$ (in yellow) are superimposed on a contour plot of $q.$
The time optimal path from $\source$ to the optimal waypoint $\sa = \argmin_{\domainT} q(\x)$ is shown in red,
with dotted red lines showing time-optimal trajectories from $\sa$ to each $\xhat_i$.  
In the left subfigure, the $\sa$-to-$\xhat_4$ trajectory overlaps the already traversed $\source$-to-$\sa$ trajectory. 
}
\label{fig:FixedT}
\end{figure}

In Figure \ref{fig:FixedT} we show two examples of such delayed-target-identification planning for $T=0.08$ and $T=0.4$ and the set of four potential targets already specified in Figure \ref{fig:determ}.  Since $q$ has multiple local minima, the value of $T$ strongly influences the optimal direction of motion starting from $\source.$  When $T$ is large enough (e.g., above $\approx 0.4528$ 
),  the optimal waypoint $\sa$ is a global minimum of $q$ on the entire $\domain.$

\begin{remark}
\label{rem:fixedT_SP}
We note that our approach in the fixed-$T$ case can be interpreted as an infinite-dimensional version of the standard Two-Stage Stochastic Programming (TSSP)  \cite{shapiro2021lectures}.
In the first stage of a {\em general TSSP} problem, some decision variable $\zeta_1 \in \G_1 \subset \R^{d_1}$ has to be chosen at a known cost $Q_1(\zeta_1)$. A random variable $\omega$ is then drawn from a known distribution and the planner chooses their second decision variable $\zeta_2 \in \R^{d_2},$ satisfying a known nonlinear constraint $\zeta_2 \in \G_2(\zeta_1, \omega)$ and incurring the cost  $Q_2(\zeta_1, \omega, \zeta_2)$.  The goal is to minimize the expected sum of these two costs; i.e., 
\[
\min\limits_{\zeta_1 \in \G_1} 
\left\{ 
Q_1(\zeta_1)  \; +  \;
\E_{\omega} \left[
\min\limits_{ \zeta_2 \in \G_2(\zeta_1, \omega)}  Q_2(\zeta_2, \omega)
\right]
\right\}.
\]

In our setting, $(\zeta_1, Q_1)$ would describe the decisions and cost prior to the target-identification time $T$, while $(\zeta_2, Q_2)$ would describe the decisions and cost after $T$.
If one needs to handle a general running cost $K$,
$\zeta_1$ becomes a choice of a path from our starting point $\x_0$, $\omega$ becomes the soon-to-be-identified target, $\zeta_2$ is a chosen path to that target
from our position $\x$ at the time $T,$ while  $q(\x)$ encodes the expected cumulative cost for the second stage. 
By using the time-dependent PDE \eqref{eq:tEik} with the terminal condition
$u = q$ on $\domain \times \{T\}$, we manage to solve this two-stage problem efficiently despite the fact that the spaces of paths $\zeta_1$ and $\zeta_2$ are infinite-dimensional.
Moreover, the value function provides the solution to {\em an infinite collection} of TSSP problems -- for all starting positions $\x_0$ simultaneously.
For a single source $\x_0$ and a constant running cost $K=1,$ one can also use the more efficient approach based on a stationary PDE \eqref{eq:Eik},  
with $\zeta_1 \in G_1 = \domainT(\source)$ interpreted as a chosen reachable waypoint and $Q_1 = T$.
\end{remark}

\vspace*{1mm}
It is also important to note that the source of initial uncertainty needs not be the random nature of the target. 
Indeed, the element of initial uncertainty could have resulted from one-time changes in the global environment -- e.g., changes in the geometry of the domain or in the running cost $K$. 
For instance, suppose that the target \change{was} known, but we had a list of possible ``obstacle'' locations (subregions of $\domain$ that would be impossible or highly risky to traverse through). If the actual obstacles were to be revealed at a later time $T$, we could use the above approach to plan the optimal path until then.
In Appendix \ref{sec:storm_fronts} we illustrate this with an example of airplane flight-path planning under weather uncertainty.


\subsection{Random discrete certainty time $T$}
\label{sec:DiscreteRandomT}

We now consider a problem in which the time of target discovery $T$ is a discretely distributed random variable. 
That is, we have a set of times $\{T_1,...,T_r\}$ with $0 < T_1 < ... < T_r,$ 
the probability $\PP(T= T_j) = p_j > 0$ is known for each $j =1, ... r,$ and 
$\sum_{j=1}^r p_j = 1.$
We will assume that $p_j$'s can be also naturally used to define conditional probabilities.
I.e., if we know that the target has not been identified by the time $T_j < T_r$, we assume that
\begin{equation}
\label{eq:cond_probab_T}
\PP \left( T = T_{j+k} \, \mid \, T > T_j \right)
\; = \;
p_{j+k} \, / 
\sum\limits_{l=j+1}^r p_l, 
\qquad \; \forall k = 1,...,r-j.
\end{equation}

In this setting, the problem can be handled by solving a sequence of HJB equations.
To avoid confusion, in this section we will always use $v$ for time-dependent value functions, saving  $u$ for solutions of static HJB PDEs.
Let $v_j(\x,t)$ denote the expected min-cost to target conditional on that target being not yet identified by the time $T_{j-1} \leq t$.  For $j=r$  this means that $\PP \left( T = T_r \mid T > t \right) = 1$,
and we are back to the case already considered in the previous section:
$v_r$ satisfies the equation \eqref{eq:tEik} on $\domain \times [T_{r-1}, T_r)$
with the terminal condition $v_r(\x,T_r) = q(\x).$  

If $j= r-1$ and $t \in [T_{r-2}, T_{r-1})$,
the target might be identified at the time $T_{r-1}$ (in which case the optimal remaining cost is encoded in $q$) or postponed until the time $T_r$ (in which case the optimal remaining cost is encoded in $v_r \left( \cdot, T_{r-1} \right)$ already found above).  Thus, on this time interval $v_{r-1}$ satisfies the same PDE \eqref{eq:tEik}, but with the terminal condition
\begin{equation}
\label{eq:time_terminal_mid}
v_{r-1} \left(\x,T_{r-1} \right) \; = \; 
\frac{p_r}{p_{r-1} + p_r} \, v_r \left(\x, T_{r-1} \right) \, + \, 
\frac{p_{r-1}}{p_{r-1} + p_r} q(\x).
\end{equation}
Iterating this process (from $j=r-1$ to $j=1$) and using $T_0= 0$ to simplify the notation, we can derive the general case.
The value function $v_j$ satisfies the PDE \eqref{eq:tEik} on
$\domain \times \left[ T_{j-1}, \, T_j \right)$
with the terminal condition
\begin{equation}
\label{eq:time_terminal_fin}
v_j(\x,T_j) \; = \; 
\left(1 -  \frac{p_j}{\sum\limits_{l=j}^r p_l} \right)  v_{j+1} \left(\x,T_j \right)
\, + \, \left(\frac{p_j}{\sum\limits_{l=j}^r p_l}\right) q(\x).
\end{equation}

As before, the computational cost of obtaining $q$ is $O(mM\log M),$ where $M$ is the number of gridpoints in the spatial grid and $m$ is the number of potential targets.
If we assume that there are $N_j$ time slices used in discretizing $[T_{j-1}^{}, T_{j}{})$, then the additional computational cost of solving the above  $r$-stage problem is
$O(M\sum_{j = 1}^r N_j).$

\begin{figure}[h]
\label{fig:T_discr}
\centering
$
\begin{array}{ccc}
\includegraphics[width = 0.285\textwidth]{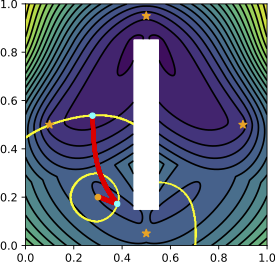} &
\includegraphics[width = 0.285\textwidth]{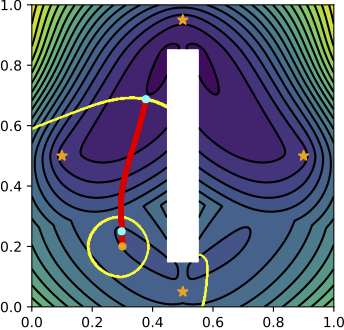} &
\includegraphics[width = 0.342\textwidth]{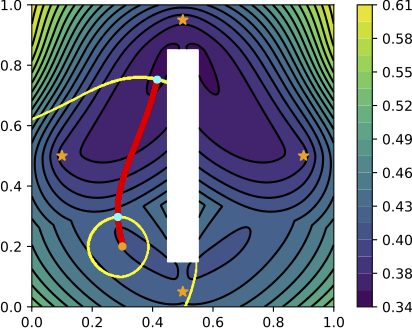}\\
(A) & (B) & (C)
\end{array}
$
\caption{Random $T$ with possible values  $T_1 = 0.08$ and $T_2 = 0.4$.  
The basic setup is the same as in Figure \ref{fig:FixedT}, but the optimal behavior is heavily dependent on $T$'s probability distribution.
\mbox{(A):  $p_1 = 0.9, \, p_2 = 0.1.$
(B):  $p_1 = 0.55, \, p_2 = 0.45.$
(C):  $p_1 = 0.1, \, p_2 = 0.9.$}
}
\end{figure}

We now consider an example using the same setup already presented in Figure \ref{fig:FixedT}, but with two possible times of target discovery: $T_1 = 0.08$ and $T_2 = 0.4$.  
Note that, even though $K=1$ and we are optimizing the expected time to target,
it is still necessary to solve two time-dependent PDEs.  As an added benefit, we recover the optimal control from every starting position, but for the sake of comparison we focus on the same
specific $\source$  (shown by an orange dot). 
In Figure \ref{fig:T_discr} the inner yellow ring illustrates $\partial \domain_{T_1} (\source),$
the boundary of the set reachable from $\source$ by the time $T_1.$ 
The cyan dot  
indicates the optimal location $\sa_1 \in \domain_{T_1} (\source)$ to reach by that time.
The outer yellow boundary shows $\partial \domain_{(T_2 - T_1)} (\sa_1),$
the boundary of the set we will have to consider if the target is not revealed at $T_1,$
and the second cyan dot $\sa_2$ is the optimal point to reach in that case.
Note that the dependence of $\domain_{(T_2 - T_1)}$  on $\sa_1$ is precisely the reason why the stationary formulation of the previous section is not usable.
Again, the background shows the level sets of $q$ (the expected time till arrival if the true target is revealed immediately). 
Figure \ref{fig:T_discr}A shows that for $p_1 = 0.9,$ the first optimal waypoint $\sa_1$ is close to $\argmin_{\x \in \domain_{T_1} (\source)} q(\x),$ the optimal solution when $p_1 = 1$ (and $T=T_1$ is fixed).   
There is a somewhat opposite effect in the case of $p_1 = 0.1$ shown in 
Figure \ref{fig:T_discr}C.  Since $T=T_2$ is far more likely, we see that $\sa_2$ is close to
$\argmin_{\x \in \domain_{T_2} (\source)} q(\x)$ and $\sa_1$ is selected (despite its relatively high $q$ value) to make this possible.
The example with $p_1 = 0.55$ presented in Figure \ref{fig:T_discr}B is more of a compromise.  It is also interesting that in this case the optimal $\sa_1$ is in the interior of $\domain_{T_1} (\source).$

The above modeling framework is flexible enough to treat possible temporal changes in the target probabilities $\bphat$ and 
even $\bphat$'s possible dependence on the certainty-time probabilities $\bm{p}$. 
We illustrate this with an additional ``emergency rescue'' example in Appendix \ref{sec:drones_example}.

\begin{remark}
We note that the discrete-random-$T$ case can be viewed as a ``multistage''(MSSP) generalization of TSSP\cite{shapiro2021lectures} mentioned in Remark \ref{rem:fixedT_SP}. 
A general $r$-stage stochastic programming problem can be written in a nested form:
\[
\min\limits_{\scalebox{0.7}{$\zeta_1 \in \G_1$}} 
\!Q_1(\zeta_1)  + 
\E_{\omega_1}  \! \!  \! \left[
\min\limits_{ \scalebox{0.7}{$\zeta_2 \in \G_2(\zeta_1, \omega_1)$} }  \!\! Q_2(\zeta_2, \omega_1) + 
\E_{\omega_2} \! \! \left[ \scalebox{0.7}{$\cdots \,$}
+ \E_{\omega_{r\text{-}1}} \! \! \! \left[
\min\limits_{  \scalebox{0.7}{$ \zeta_{r} \in \G_{r}(\zeta_{r\text{-}1}, \omega_{r\text{-}1})$} }
\!\!\!\! Q_{r}(\zeta_{r}, \omega_{r\text{-}1})
\right]
\scalebox{0.7}{$\cdots$} 
\right]
\right] \! ,
\]
with $\omega_1,\omega_2,\cdots,\omega_{r-1}$ being a stochastic process.
The expectation $\E_{\omega_j}[\cdot]$ is conditioned on an observed realization of $\omega_1,\omega_2,\cdots,\omega_{j - 1}$.

For our model, the random outcomes of whether the true target is revealed at each $T_j$ correspond to the stochastic process of MSSP.
The terminal conditions defined in \eqref{eq:time_terminal_mid} and \eqref{eq:time_terminal_fin} are conditional expectations of the two possibilities at $T_j$, provided that $T > T_{j - 1}$.
The HJB solution $v_j(\x, t)$ in the $j$th stage encodes the optimal remaining cost for any starting state $\x$, which \change{is} naturally used to determine the terminal conditions for the previous stage.  Since we are obtaining the optimal control for all initial $(\x,t)$ simultaneously, this is again equivalent to solving an infinite family of infinite-dimensional MSSP problems.
\end{remark}


\subsection{Exponentially distributed certainty time $T$}
\label{sec:ExponentialRandomT}

We now move beyond problems that can be solved stage-by-stage in a typical MSSP fashion.
Section \ref{sec:DiscreteRandomT} assumes the certainty time to be a discretely distributed random variable.
A natural extension is to consider a continuously distributed $T$.  
The total cost can be again split into two parts (before and after certainty is achieved), with the expected cost of the latter part
still encoded by $q(\x)$ defined in section \ref{sec:FixedT}.
The target identification is the termination of the first part, and the optimal control of such {\em randomly-terminated} processes has already been
considered in  \cite{andrews2013deterministic}.  We will follow the same approach here, taking $T$ to be exponentially distributed with a known rate $\lambda >0.$  I.e.,   $\E(T) = 1/\lambda$ and $\PP(T > t + \tau \, \mid \, T \geq t) = e^{-\lambda \tau}$ for all $t,\tau >0.$

Under this assumption, the expected cost can be computed as
\begin{align*}
\J(\x, \ba(\cdot)) \, & = \, \int_0^{\infty} \lambda e^{-\lambda T} \left[  \int_0^T K(\y(t))  \, dt \, + \, q\left(\y(T)\right) \right] dT \\
& = \int_0^{\infty} e^{-\lambda t}\left[ K(\y(t)) + \lambda q(\y(t)) \right] \, dt.
\end{align*}
Thus, this problem can be also viewed as a discounted infinite-horizon problem with $\lambda$ interpreted as a discounted factor and $\left( K(\x) + \lambda q(\x) \right)$ as a new running cost.
The value function $u^{\lambda}(\x) = \inf_{\ba(\cdot)}\J(\x, \ba(\cdot))$ can be recovered as a \textit{viscosity solution} of a Hamilton-Jacobi PDE 
\begin{equation}
\label{eq:random_terminate}  
\lambda \left( u^{\lambda}(\x) - q(\x) \right)  \, + \, |\nabla u^{\lambda}(\x)| f(\x) \; = \; K(\x).
\end{equation}

Focusing on time optimality, we can use $K=1$ to have $u^{\lambda}(\x)$ encode the (minimized) full  expected-time-to-target from $\x$ to uncertain $\xhat$.  Alternatively, we can use the fact that  $\E(T) = 1/\lambda$ is trajectory independent and take $K=0$ to have $u^{\lambda}(\x)$ encode the (minimized) expected-time-to-target-after-it-is-identified.  We follow \cite{andrews2013deterministic} and use the latter, which 
yields $u^{\lambda} \leq q$ since one of the options is to stay in place until the target is identified.
\change{
This allows recovering the value function from a quasi-variational inequality of obstacle type
\begin{equation}
\label{eq:random_terminate_qvc}
\max \left(
u^{\lambda} - q,
\quad
\lambda \left( u^{\lambda} - q \right)  +  f |\nabla u^{\lambda}| 
\right) \; = \; 0.  
\end{equation}
}
We use $\M^{\lambda}$ to denote the {\em motionless set} of starting positions at which $u^{\lambda} = q,$
\change{so that} 
$\lambda \left( u^{\lambda} - q \right)  +  f |\nabla u^{\lambda}| = 0$ is solved on $\domain \backslash \M^{\lambda}.$
For any starting position $\source \not \in \M^{\lambda},$ the optimal trajectory is found by gradient descent in $u^{\lambda}$ and leads to some motionless waypoint $\sa \in \M^{\lambda}.$  If the target is not identified by the time we reach $\sa,$ it is optimal to stay there until $\xhat$ is revealed.  
This is a problem with {\em free boundary} since $\M^{\lambda}$ is not known in advance.  (With $K=0$, it is easy to show that every point in $\M^{\lambda}$ is a local minimum of $q.$  The converse is true for global minima of $q$, but its local minima may or may not be motionless depending on $\lambda$ and the global properties of $q.$)  Nevertheless, a modified version of Fast Marching Method can be still used to compute the numerical approximation efficiently (see \cite[Section 3.3 and Appendix B]{andrews2013deterministic} for implementation details).
Since $q$ has to be computed first before solving  \eqref{eq:random_terminate}, the overall cost is  $O\left( (m+1)M\log M \right).$ 

\begin{figure}[h]
\centering
$
\begin{array}{ccc}
\includegraphics[width = 0.305\textwidth]{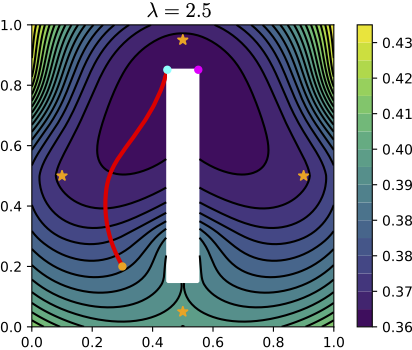}&
\includegraphics[width = 0.305\textwidth]{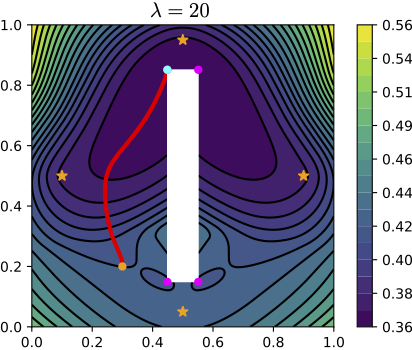}&
\includegraphics[width = 0.305\textwidth]{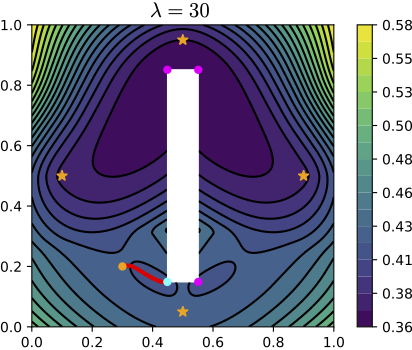}\\
(A) & (B) & (C)
\end{array}
$
\caption{Exponentially distributed $T$.
Each subfigure shows the level sets of $u^{\lambda},$ starting position $\source$ (in orange), optimal waiting position $\sa$ (in cyan), motionless local minima of $u^\lambda$ (in magenta) and the optimal trajectory (in red). 
(A): with $\lambda = 2.5$ an early target identification is not very likely, and the optimal trajectory leads to the global minimum of $q$. 
(B): with $\lambda = 20$ still heading to the global minimum although the nearby local minimum is already motionless.
(C): with $\lambda = 30$ reaching the global minimum before the target identification is less unlikely, and it is optimal to head toward the 
local minimum nearby.}
\label{fig:exponential}
\end{figure}

The value of parameter $\lambda$ determines
whether one should make a quick transition to a nearby local minimum or spend more time on a path to the global minimum 
(risking that the target is identified while we are still traveling through high values of $q$).
The example in Figure \ref{fig:exponential} illustrates this, with $\E[T] = 1/2.5 = 0.4$ in the left subfigure producing 
a trajectory similar to the one in Figure \ref{fig:FixedT}B, and  $\E[T] = 1/30$ in the right subfigure resulting in a completely different path to a local minimum.
We note that with an intermediate $\E[T] =  1/20$ that nearby local minimum is already motionless, but the optimal trajectory from our chosen $\source$ still leads to a much farther global minimum. 



\section{Robust path planning}
\label{sec:Robustness}

We now turn our attention to robust control under initial uncertainty, using the same example previously considered in section \ref{sec:FixedT} (i.e., minimizing time-to-target that will be revealed at a known time $T$). 
Our goal is to compare several different notions of robustness: the worst-case optimization (\S \ref{ss:worst_case}), 
the risk-sensitive average optimization (\S \ref{ss:risk_sensitive}), 
optimizing the average with a hard constraint on the worst case (\S \ref{ss:average_worst}),  
and optimizing the average with a ``chance constraint'' on bad outcomes (\S \ref{ss:probabilistic}).
As we will show,
the corresponding optimal waypoints and trajectories are quite different in each case, including the probabilistic definition of optimal waypoints in \S \ref{ss:probabilistic}.
Along the way, we will also provide a brief discussion of how these robustness models relate to each other.
Two of these (in \S  \ref{ss:worst_case} and \S \ref{ss:average_worst}) will be also extended to the case of an exponentially distributed random certainty time $T$ (introduced in \S \ref{sec:ExponentialRandomT}).

All of the above describe robustness to {\em aleatoric uncertainty}; i.e., uncertainty due to a random choice based on a known probability distribution $\bphat$ over a known
set of potential targets $\Xhat.$
We end this section by exploring robustness with respect to  {\em epistemic uncertainty} (\S \ref{ss:epistemic}); i.e., possible uncertainties in $\bphat$ or $\Xhat.$


\subsection{Worst\change{-}case scenario}
\label{ss:worst_case} 

Perhaps the most obvious robust approach is to consider a Stackelberg-type ``game against nature'', where the worst target is always chosen  (by a supposed opponent) based on our choice of a waypoint.  In this case, the time to target after identification is
\begin{equation}
\qbar(\x) = \max\{u_1(\x),\cdots,u_m(\x)\}. 
\end{equation}
and the worst-case optimal solution is to go from $\source$ to any waypoint  
\begin{equation}
\sw \in \argmin\limits_{
\x \in \domainT(\source) }  \qbar(\x).
\end{equation}

\begin{figure}[h]
\centering
$
\begin{array}{cc}
\includegraphics[width = 0.44\textwidth]{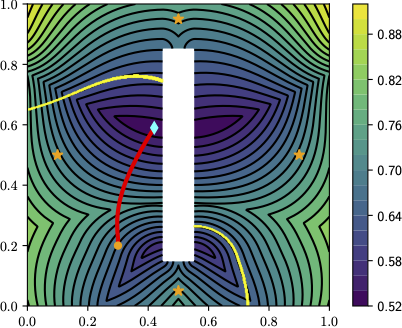} & \quad
\includegraphics[width = 0.442\textwidth]{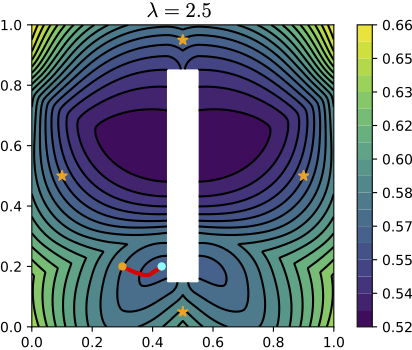}\\
(A) & (B)
\end{array}
$
\caption{Worst-case optimal planning.
(A): the same (``Fixed $T=0.4$'') setting as in Figure \ref{fig:FixedT}B. 
The cyan diamond represents a $\sw$ and the red curve is the time-optimal trajectory, with level sets of $\qbar$ in the background.
(B): the same (``Exponentially distributed $T$ with $\E[T] = 0.4$'') setting as in Figure \ref{fig:exponential}A except that $q$ is replaced by $\qbar$. Level sets of $u^\lambda$ are shown in the background.
} 
\label{f:worst_case}
\end{figure}

This approach (illustrated in Figure \ref{f:worst_case}A) does not use any information from $\bphat$ and is thus very conservative.  E.g., the worst\change{-}case scenario could be very unlikely and $q(\sw)$ can be quite far from the average-case optimal $q(\sa).$  Similarly, the worst\change{-}case scenario $\qbar(\sa)$
can be much costlier than $\qbar(\sw).$  While these quantities can be compared directly, a more nuanced discussion of the ``worst versus average'' trade-offs will be covered in \S \ref{ss:average_worst}.
The computational cost of obtaining $\sw$ through our Fast Marching implementation is still $O(mM\log M)$.

As we show in Figure \ref{f:worst_case}B, the same approach also works for the exponentially distributed random certainty time discussed in section \ref{sec:ExponentialRandomT}.   For the specified rate $\lambda > 0$, the worst-case optimal trajectory is found by using $\qbar$ instead of $q$ in PDE \eqref{eq:random_terminate} and then following $(- \nabla u^\lambda)$ from $\source$ to the corresponding motionless point.  This waypoint $\sw$ is one of the local minima of $\qbar,$ but we may not reach it if the target is identified earlier.  (Note that this is the worst case with respect to $\xhat$ only while $T$ is still assumed to be random.
If we allowed $T$ to be chosen by the opponent/nature,  it would always result in $T=0,$ yielding $\qbar(\source)$ as our time to the worst\change{-}case target.)

\subsection{Risk-sensitive optimization}
\label{ss:risk_sensitive} 

\begin{figure}[h]
\centering
$
\begin{array}{cc}
\includegraphics[width = 0.45\textwidth,valign=t]{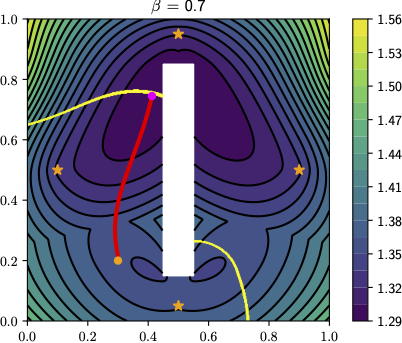} & \quad
\includegraphics[width = 0.45\textwidth,valign=t]
{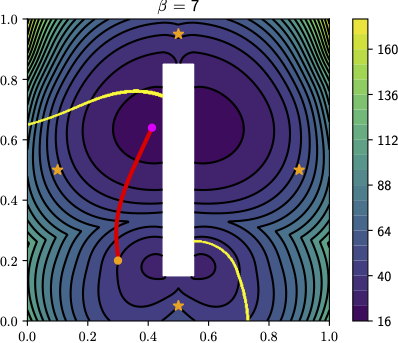}\\
(A) & (B)
\end{array}
$
\caption{Risk-sensitive planning. The red curve is a time-optimal trajectory toward $\srs$ (magenta dot). Level sets of $q_\beta$ are shown in the background. 
The position of $\srs$ depends on the magnitude of $\beta$. 
(A): A small $\beta$ leads to $\srs$ being close to $\sa$ (compare with Figure \ref{fig:FixedT}B).
(B): $\beta$ is large and $\srs$ moves closer to $\sw$ (compare with Figure \ref{f:worst_case}A).
}
\label{fig:risk_sensitive}
\end{figure}

The notion of risk sensitivity most commonly used in optimal control and Markov Decision Processes is based on minimizing the expectation of an exponential function of the outcome \cite{HowardMatheson_1972, FlemingSoner}.
Given a risk-sensitivity parameter $\beta > 0,$ it is natural to define
\begin{equation}
\qbeta(\x) = \E_{\xhat} [e^{\beta u(\x)}] = \sum_{i=1}^m \phat_i e^{\beta u_i(\x)}.
\end{equation}
The waypoint is then chosen as 
\[
\srs \in \argmin\limits_{
\x \in \domainT }  \qbeta(\x). 
\]
Figure \ref{fig:risk_sensitive} illustrates this approach for $T=0.4$ and two different $\beta$ values. 
Higher $\beta$ corresponds to more risk-averse behavior. 
As $\beta$ increases, the exponential of the worst outcome starts to dominate the expectation; so, $\srs \rightarrow \sw$ as  $\beta \rightarrow +\infty.$  On the other hand, when $\beta$ is small, the Taylor expansion yields $\qbeta \approx 1 + \beta q$; so, $\srs \rightarrow \sa$ as  
$\beta \rightarrow 0.$   

Another approximation valid for small $\beta$'s establishes the connection with the mean-variance trade-offs:
\begin{equation}
\label{eq:mean_variance}
\log (\qbeta ) / \beta \; \approx \; \E_{\xhat} u(\x) \, + \, 
\frac{\beta}{2} \E_{\xhat} \left[ u(\x)   -  \E_{\xhat} u(\x) \right]^2.  
\end{equation}
Unfortunately, this holds only asymptotically and is not very useful in estimating the worst-average case trade-offs for positive $\beta$ values; see Figure \ref{f:worst_expected} and Remark \ref{rm:jensen}.

\subsection{Optimizing the average with a hard constraint on the worst case}
\label{ss:average_worst} 

\begin{figure}[h]
\centering
$
\begin{array}{cc}
\includegraphics[width = 0.45\textwidth]{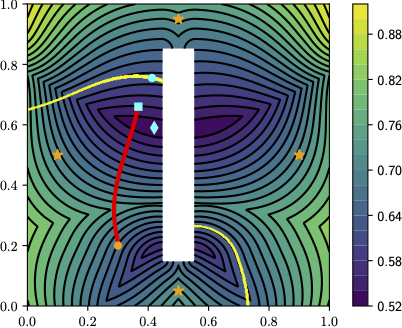} & \quad
\includegraphics[width = 0.47\textwidth]{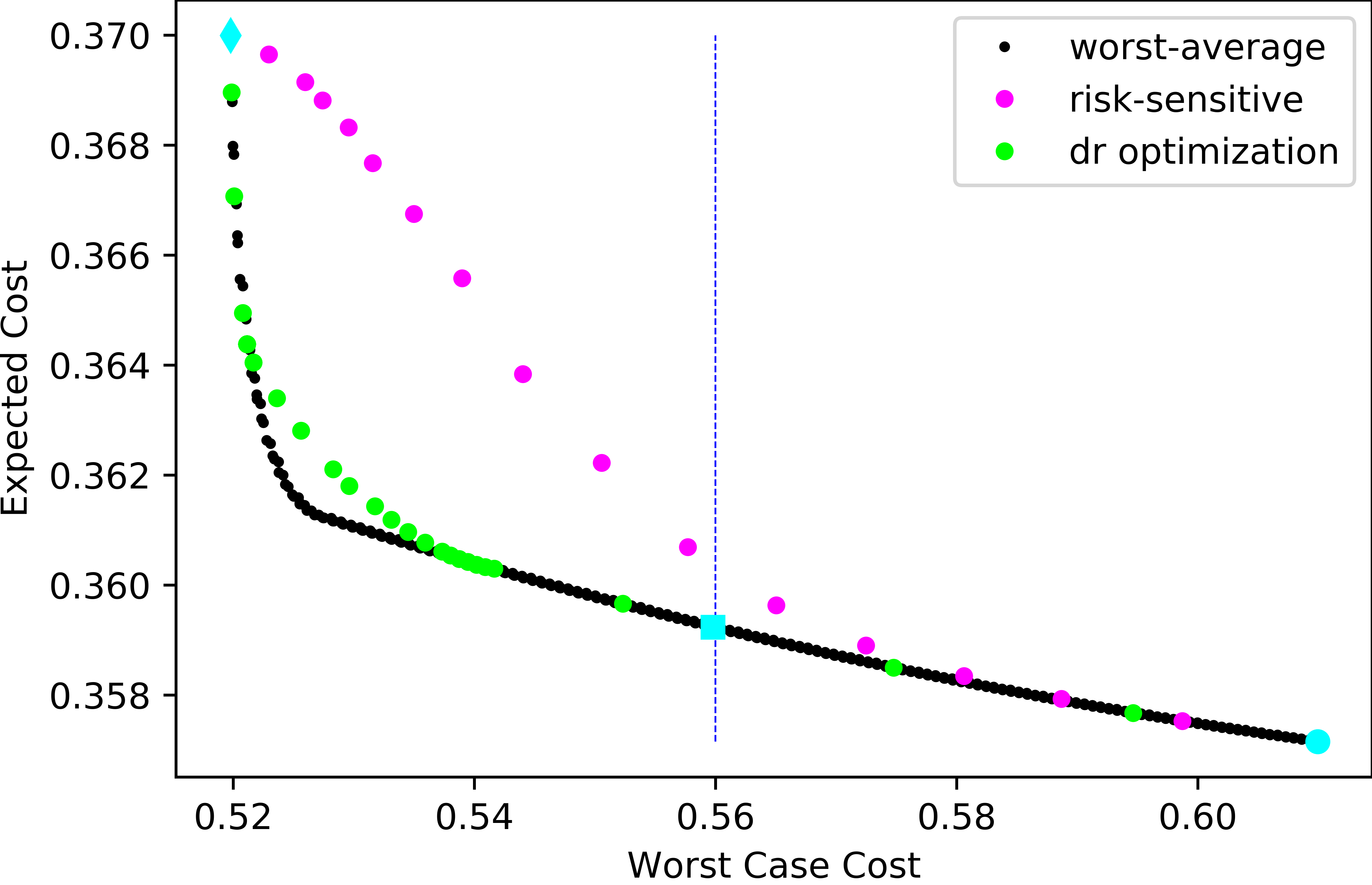}\\
(A) & (B)
\end{array}
$
\caption{Minimizing averaged cost with a hard constraint $\worst = 0.56$. 
(A):
Three different optimal waypoints ($\sa,  \sw, \saw$) corresponding to sections \ref{sec:FixedT}, \ref{ss:worst_case} and \ref{ss:average_worst} are shown by cyan markers (the circle, diamond and square respectively).
The red curve is the time optimal trajectory from $\source$ to the latter, with the level sets of $\qbar$ in the background.
(B): The Pareto Front approximation is shown in black, with the cyan markers corresponding to the 3 waypoints shown on the left.
The dashed line is $\qbar = 0.56$ and the square marker is the lowest one on \change{the} Pareto front satisfying $\qbar \leq 0.56$.
The scattered magenta and green dots correspond to the risk-sensitive and distributionally-robust waypoints discussed in sections
\ref{ss:risk_sensitive} and \ref{ss:epistemic}.
}
\label{f:worst_expected}
\end{figure}

As usual\change{,} when evaluating more than one performance criterion, we have to employ the notion of {\em Pareto-optimality}.
Focusing on $K=1$ and deterministic $T$,
we will say that a waypoint $\z_1$ is dominated by another waypoint $\z_2$ if both $\qbar(\z_2) \leq \qbar(\z_1), \; q(\z_2) \leq q(\z_1),$
and at least one of these inequalities is strict.  We will say that a waypoint $\z \in \domainT(\source)$ is Pareto-optimal if it is not dominated by any other  waypoint in $\domainT(\source).$  Collectively, these waypoints encode all rational ``average versus worst'' trade-offs since none of them can be improved with regard to both criteria simultaneously.  These trade-offs are well represented visually by a Pareto Front 
$PF = \left\{ \left(\qbar(\z), q(\z) \right) \, \mid \, \z \text{ is Pareto-optimal} \right\},$
which a practitioner could use to find the desired compromise.
Every such waypoint can be viewed as an answer to a constrained optimization problem: minimize $q$ subject to $\qbar \leq \worst$ (or, equivalently,
ensuring that the total time to target is below $T + \worst$).
For a fixed starting position $\source$, this can be done by restricting the domain to
\[
\omegatc(\source) \, = \, \omegat \bigcap \omegac,    \qquad \text{where } \omegac \, = \, \{ \x \in \Omega \, \mid \,  \qbar(\x) \leq \worst\},
\] 
and then selecting a waypoint 
\[
\saw \in \argmin\limits_{
\x \in \omegatc(\source) }  q(\x).
\]
Note that the problem has no solution if $\worst$ is too small (resulting in $\omegatc(\source) = \emptyset$), while for large enough $\worst$
the problem becomes unconstrained, yielding the same waypoints obtained in section \ref{sec:FixedT}.
A more efficient approximation of the entire $PF$ is obtained by evaluating $(\qbar, q)$ at every gridpoint in  $\omegat(\source)$
and then removing the dominated ones.  The result of this approach for $T=0.4$ and $C=0.56$ is shown in Figure \ref{f:worst_expected}.

\begin{remark}
\label{rm:jensen}
By the optimality of $\srs$ and Jensen's inequality, we have
\[
\E_{\xhat} [e^{\beta u(\bs)}] = \qbeta(\bs) \geq \qbeta(\srs) =  \E_{\xhat} [e^{\beta u(\srs)}] \geq e^{\beta \E_{\xhat} [u(\srs)]} = e^{\beta q(\srs)}.
\]
After rearranging this inequality and applying \eqref{eq:mean_variance} we obtain
\[
q(\srs) \leq \log\left( \E_{\xhat} [e^{\beta u(\bs)}] \right)/\beta \approx q(\bs) + 
\frac{\beta}{2} \E_{\xhat} \left[ u(\bs)   -  q(\bs) \right]^2.
\]
which means the vertical gap between $q(\srs)$ and $q(\bs)$ is approximately bounded by $\frac{\beta}{2} \text{Var}[ u(\bs) ]$ when $\beta$ is relatively small.
\end{remark}

\begin{figure}[h]
\centering
$
\begin{array}{ccc}
\includegraphics[width = 0.305\textwidth]{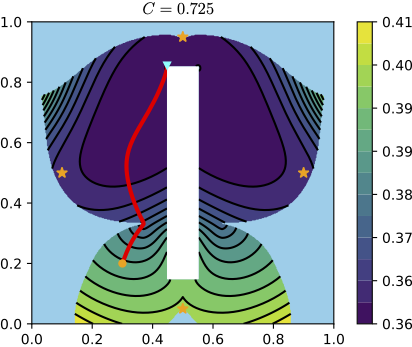} &
\includegraphics[width = 0.305\textwidth]{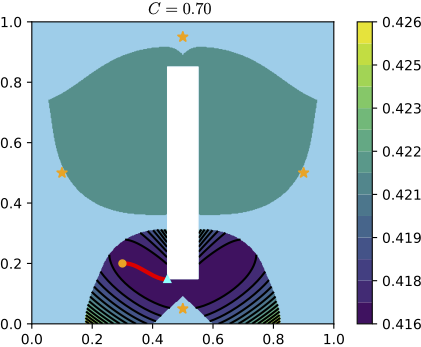} &
\includegraphics[width = 0.31\textwidth]{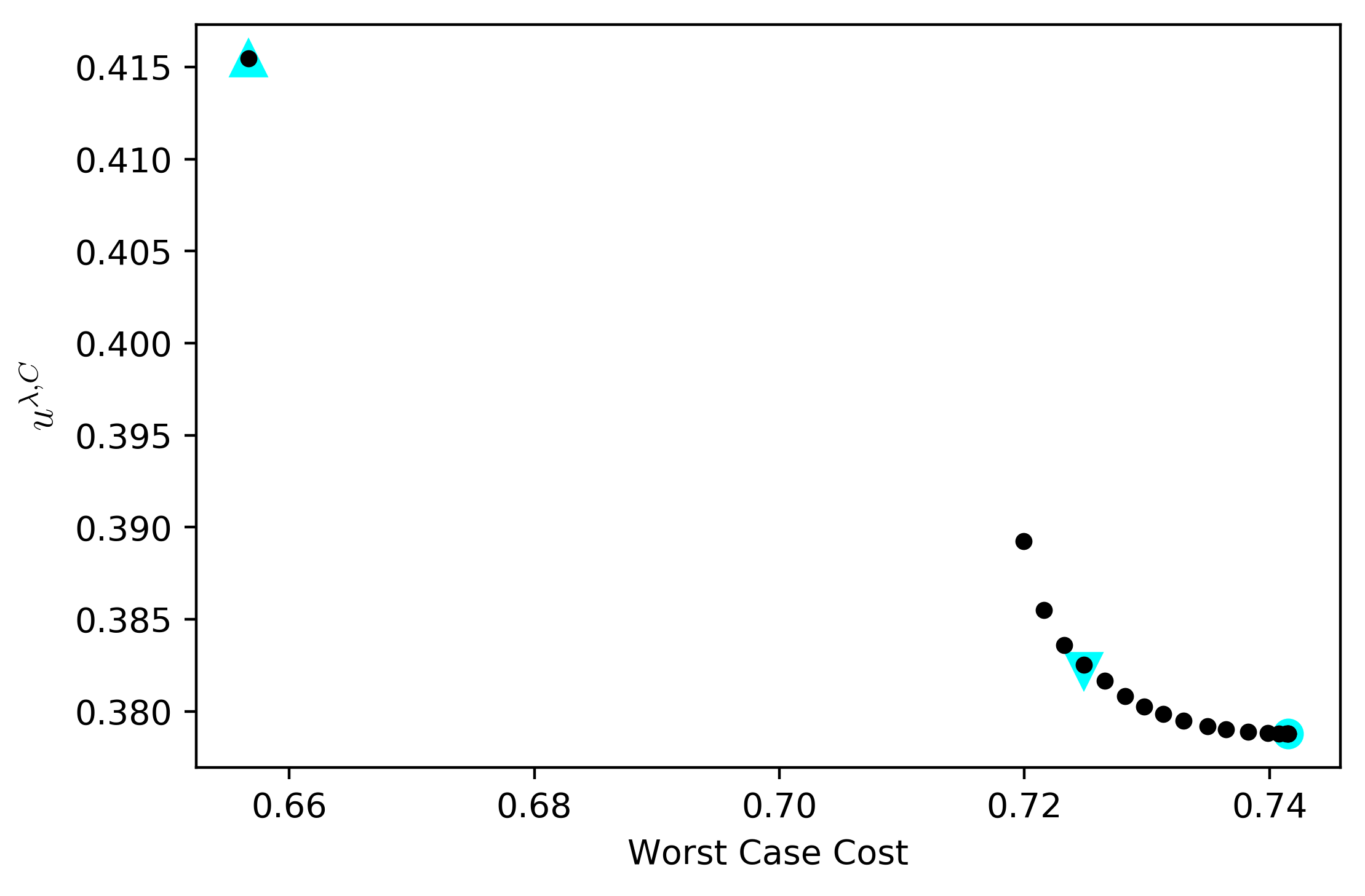} \\
(A) & (B) & (C)
\end{array}
$
\caption{Worst-case-constrained optimiziation for the exponentially distributed $T$ with $\lambda = 2.5$.
(A): a binding constraint $\qbar \leq 0.725.$ 
The cyan dot shows the same waypoint as in Figure \ref{fig:exponential}A; i.e., $\saw = \sa$ is the global minimum of $q$.  But the optimal trajectory bends to stay inside $\omegac$ and $u^{\lambda,C} (\sa) > u^{\lambda} (\sa).$
(B): a non-binding constraint $\qbar \leq 0.7.$  As $C$ decreases, $\omegac$ becomes disconnected and the optimal trajectory leads to a closer local minimum.
The dark cyan region is another connected component of $\omegac$. 
(C): the Pareto Front (worst vs average case performance). 
There are no solutions for $C < \qbar(\bs) \approx 0.6566$  
and we use $C$ values between $0.66$ and $0.81$ to approximate PF.
The cyan upward and downward triangles correspond to the 
the left and middle figures respectively. 
The cyan dot corresponds to the unconstrained $\lambda$-optimal trajectory in Figure \ref{fig:exponential}A. 
}
\label{f:exponential_hardconst}
\end{figure}

If $\worst$ and $T$ are fixed, the solution can be obtained for all starting positions at once by solving the time-dependent PDE \eqref{eq:tEik} with $K=1$ and the terminal condition  $u(\x,T) = q(\x)$ on $\omegac$ and $u(\x,T) = +\infty$ on $\domain \backslash \omegac.$  
A constraint on the {\em remaining time-to-target after $\xhat$ is identified} can be similarly handled even with the exponentially distributed random $T$  discussed in section \ref{sec:ExponentialRandomT}.  We define $u^{\lambda,C}$ to be the solution of \eqref{eq:random_terminate} on the restricted domain $\omegac$  (see Figures \ref{f:exponential_hardconst}A, \ref{f:exponential_hardconst}B for two examples with different $\worst$ values).  
By varying $C$ we can also approximate the Pareto Front as shown in Figure \ref{f:exponential_hardconst}C.  
The average case performance is encoded by $u^{\lambda,C}(\source)$, while the worst case corresponds to the 
largest value of $\qbar$ encountered on the optimal trajectory from $\source$ (found by following  $-\nabla  u^{\lambda,C}$).
This worst\change{-}case result is equal to $\worst$ if that trajectory touches $\partial \omegac \backslash \boundary,$
but as Figure \ref{f:exponential_hardconst}B shows, this is not always the case.

Before moving on to other types of robust path planning, we note that this ``constrain the worst, optimize the average'' strategy can be also computationally efficient even in a more general setting, where the uncertainty is not monotonically decreasing.  E.g., it was successfully used in robust routing on stochastic networks \cite{ermon2012probabilistic}, where 
the target was pre-specified, but each edge transition incurred a random penalty.

\begin{remark}
The computational cost of producing $q, \qbar,$ and $\qbeta$ is the same.
So, finding $\sw$ or $\srs$ for a specific $\beta$ costs exactly as much as obtaining the risk-neutral $\sa.$
The same is also true for $\saw$ corresponding to any specific $C$.  (In fact, the computational cost is usually reduced since  $\omegac$ is smaller than $\Omega.$)  
To approximate the Pareto front, we can first sort all $(\qbar, q)$ pairs in $\omegat(\source)$
by the $\qbar$ values, then go through all of them in increasing order 
and prune all dominated pairs.  The computational cost of this process is $O(\tilde{M} \log \tilde{M}),$
where $\tilde{M} \leq M$ is the number of gridpoints 
in $\omegat(\source)$.
\end{remark}


\subsection{Probabilistic constraints on bad outcomes}
\label{ss:probabilistic} 
Probabilistic constraints (also frequently known as {\em chance constraints}) provide a fairly general version of robustness when optimizing the expected performance.  
Unfortunately, they often result in much tougher optimization problems \cite{white1974dynamic, pfeiffer2018two, palmer2017optimal}, but we will show that in our simple setting the additional computational cost is  
only moderate.

Focusing on each potential target $\xhat_i$, it is natural to define its ``$\worst$-unreachable set''
$\hat{\domain}_i = \{\x\in \omegat(\source) \,  \mid \,  u_i(\x) > \worst \}$ and 
the corresponding indicator function $\chi_i^{}(\x) = \mathbf{1}_{\hat{\domain}_i}(\x)$.
The risk function can be then defined as the probability of not meeting this constraint on the remaining time to target if we are at the point $\x$ when that target is finally identified.  I.e.,
\begin{equation}
\label{eq:risk}
r(\x) = \sum_{i=1}^m \phat_i \chi_i^{} (\x).
\end{equation}
We note that the worst-case-constrained expectation optimization of section \ref{ss:average_worst} can be viewed 
as minimizing $q(\x)$ while guaranteeing that $r(\x) = 0$.  This approach is highly conservative since it disqualifies many waypoints $\x$ with small $q$ simply because some $u_i (\x) > \worst$ even if the corresponding $\phat_i$ is quite small.
It might seem natural to relax this hard constraint by minimizing $q$ over the set $\{\x \in \omegat(\source) \,  \mid \,  r(\x) \leq \epsilon\}$.
But since our goal is to limit the overall risk, better results are obtained by allowing the planner to use probabilistic/mixed strategies.  
We can occasionally use waypoints with $r > \epsilon$ which are attractive because of their small $q$, and still limit the overall risk by also using waypoints with $r < \epsilon$ whose $q$ values are not as good.  
An intermittent use of both $\sa$ and $\saw$ is a simple example of this kind of strategy.  With $\worst = 0.56$ and $T= 0.4$ corresponding to Figure \ref{f:worst_expected}, this actually happens to be optimal if we use  $\sa$ with probability $\theta_{\sa} = \min (1, \epsilon / r(\sa) )$ and head to 
$\saw$ with probability $(1- \theta_{\sa})$.  However, we will show that this use of $\sa$ and $\saw$ is not generically optimal.  In many cases,
the set $\omegatc(\source)$ may be even empty, making  $\saw$ undefined.  Below we consider a more interesting example with a lower $\worst$ value resulting in a positive risk $r$ on the entire $\omegat(\source).$

For full generality, we could let a planner select a probability measure over all possible waypoints in $\omegat(\source)$,
minimizing the expected $q$ and restricting the expected $r$ with respect to that measure.  But in practice we will assume that
the planner selects a discrete pdf $\bm{\theta}=(\theta_1, \ldots, \theta_n)$ over the set of gridpoints $\mathcal{X} = \{\x_1,\cdots,\x_n\} \subset \omegat(\source).$
I.e., we will assume that $\theta_j \geq 0$ is the probability of choosing a gridpoint $\x_j$ as our waypoint.
In this setting, the probabilistically-constrained optimization can be performed by solving a finite-dimensional {\em linear program} with
$\theta_j$'s as decision variables:
\begin{equation}
\label{eq:lp_probab_constr}
\begin{array}{ll@{}ll}
\text{minimize}  & \displaystyle\sum\limits_{j=1}^n \theta_j  q(\x_j) &\\
\text{subject to}& \displaystyle\sum\limits_{j=1}^n    r(\x_j) \theta_j \leq \epsilon, &\\
                       &  \displaystyle\sum\limits_{j=1}^n \theta_j = 1, \quad \theta_j \geq 0,\ j = 1,\cdots,n. 
\end{array}
\end{equation}
As long as the feasible set is non-empty, it will be an $(n-1)$-dimensional polytope with a minimum attained at one of its vertices (i.e., where $(n-1)$ constraints are active).  Since there is only one constraint based on the risk, at least $(n-2)$  values of $\theta_j$'s are actually zero at any vertex.  Thus, regardless of $m$ and $n$, there always exists an optimal solution assigning non-zero probability to at most 2 waypoints as long as the feasible set is non-empty.  If one wished to impose $\kappa$ probabilistic constraints
\[
\sum\limits_{j=1}^n    \PP(u > \worst_k) \leq \epsilon_k,
\qquad
k=1, ... \kappa,
\]
a similar argument shows that an optimal mixed strategy could be found based on at most $(\kappa+1)$ waypoints\footnote{The authors are grateful to Dmitriy Drusvyatskiy for pointing this out.}.

\begin{figure}[h]
\centering
$
\begin{array}{cc}
\includegraphics[width = 0.46\textwidth]{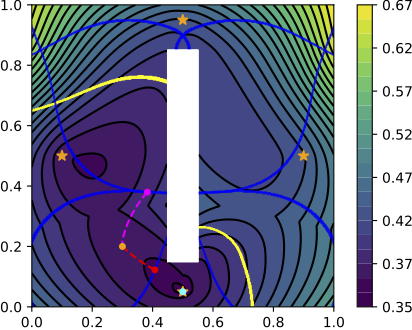} & \quad
\includegraphics[width = 0.45\textwidth]{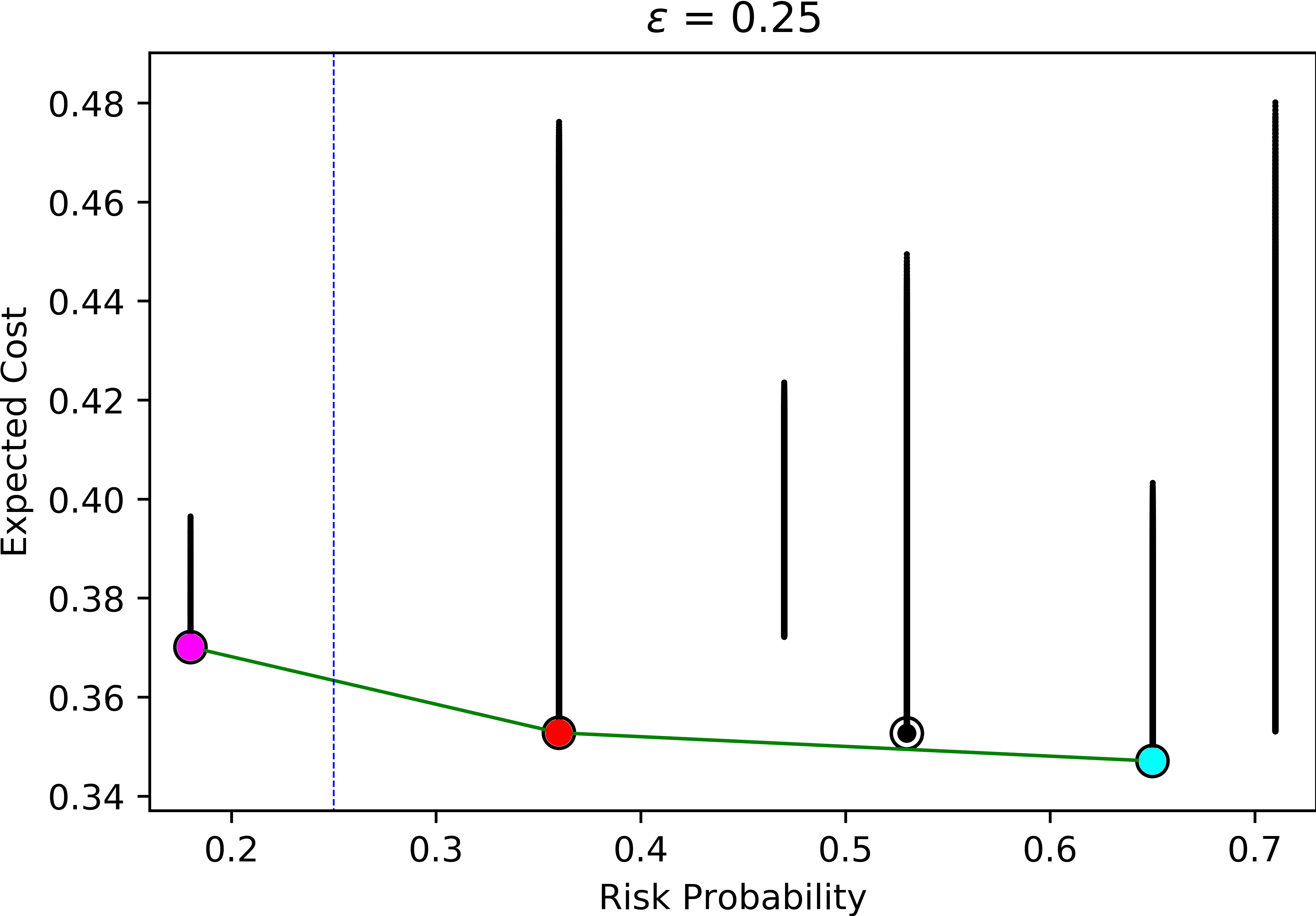}\\
(A) & (B)
\end{array}
$
\caption{Probabilistic constraints ($C = 0.365; \; \epsilon = 0.25)$. 
(A): Level sets of $q(\x)$ corresponding to $\phat_1 = \phat_2 = 0.18, \phat_3 = 0.35,$ and $\phat_4 = 0.29$
with the same target locations as in Figure \ref{fig:FixedT}.  The boundary of $\domain_T(\source)$ is yellow and
the boundaries of $\hat{\domain}_i$'s are blue.
The ``unconstrained optimal'' waypoint 
$\bs = \argmin_{\x \in \domain_T(\source)} q(\x)$
 is located exactly at the target $\xhat_3$ (in cyan).
The pair of $(C,\epsilon)$-optimal waypoints are shown in red and magenta,
and the optimal mixed strategy selects them with the respective probabilities of $0.3889$ and $0.6111$.
(B): The set $\mathcal{X}^{r,q}$ \change{is shown} in black, its non-dominated points \change{are} shown by circular markers, while $\V$ (the ``southwestern boundary'' of its convex hull $\Upsilon$) is shown in green.  The line $r=\epsilon$ is shown in blue and the ratio in which it divides the segment of $\V$ determines the probabilities of red and magenta waypoints in the optimal mixed strategy.
The magenta waypoint can be also viewed as maximizing the probability of a desirable event (i.e., $u < C$).
Even though $m=4$, there are only 6 values of $r$ in this example since many of the $2^m$ groups of gridpoints are empty
and a few groups give the same values (since $\phat_1 = \phat_2$).
}
\label{fig:prob_constraint}
\end{figure}

While the general linear programming algorithms are certainly suitable for solving \eqref{eq:lp_probab_constr} based on a 2D grid, we also describe a better (more geometric and efficient) approach based on $(r,q)$-Pareto optimality.
The set of possible waypoints $\mathcal{X}$ is naturally mapped to 
$\mathcal{X}^{r,q}  = \left\{ \left(r(\x), q(\x)\right) \, | \, \x \in \mathcal{X} \right\}.$
Any point within its convex hull $\Upsilon = co \left(\mathcal{X}^{r,q} \right)$ is attainable by selecting an appropriate $\bm{\theta}$.
The solution of \eqref{eq:lp_probab_constr} corresponds to 
$ \min \left\{ q \, | \, (r,q) \in \Upsilon, \, r \leq \epsilon \right\},$
which is actually attained on $\Upsilon$'s ``southwestern'' boundary 
\[
\V = 
\left\{ (r,q) \in \partial \Upsilon \, | \,  
\text{ with } r < r' \text{ or } q < q', \,\forall (r',q') \in  \Upsilon \setminus \{ (r,q)\} \right\}.
\] 
This $\V$ is a simple polygonal chain whose vertices are among the Pareto-optimal points of $\mathcal{X}^{r,q},$
stretching from a ``first-$r$-then-$q$'' minimizing waypoint on the left to a ``first-$q$-then-$r$'' minimizing waypoint on the right; see Figure \ref{fig:prob_constraint}B.  
Any of the planar convex hull methods can be used to compute $\V$ efficiently; e.g., \textit{Graham scan} \cite{graham1972efficient} or \textit{monotone chain} \cite{andrew1979another} algorithms can compute it in $\mathcal{O}(n \log n)$ operations.  If 
$\partial \Upsilon$ 
has only $k \ll n$
vertices on it, \textit{Chan's algorithm} \cite{chan1996optimal} yields an even better worst\change{-}case complexity of $\mathcal{O}(n \log k).$

Once $\V$ is obtained as an $r$-sorted list of its $k$ vertices, this essentially solves \eqref{eq:lp_probab_constr} for all possible $\epsilon$'s at once.  If $\epsilon < \min\limits_{(r,q)\in\V} r$, the problem has no solution.  
If $\epsilon \geq \max\limits_{(r,q)\in\V} r,$  the problem is unconstrained.  
In all other cases, the line $r=\epsilon$ intersects some segment of $\V$, which can be found by a binary search in 
$\mathcal{O}( \log k)$ steps.
This intersection point yields the corresponding probabilities with which the waypoints at the segment vertices should be selected; see Figure \ref{fig:prob_constraint}B.  If the line $r=\epsilon$ passes through a vertex of $\V$, it is optimal to use the corresponding single waypoint deterministically.

\begin{remark}  
For the case $m \ll \log n,$ additional speed up can be achieved by ignoring the dominated (i.e., the non-Pareto-optimal) points in $\mathcal{X}^{r,q}$.    
Note that in our setting $r(\x)$ is a piecewise constant function with at most $2^m$ possible values
(based on whether or not $\x$ belongs to each $\hat{\domain}_i$).  
When computing $r(\x)$, we can check which $\hat{\domain}_i$'s $\x$ belongs to and compare its $q(\x)$ with the smallest $q$ already computed within the same-$r$-value group.
Thus, maintaining a list of minimal $q$ values for all groups takes $\mathcal{O}(nm)$ steps.
Notice that all the $2^m$ $r$ values are generated by taking sums of all the subsets of $\{\phat_1,\cdots,\phat_m\}$.
If all the $\phat_i$'s are sorted (which can be done with $\mathcal{O}(m\log m)$ steps), using \change{the} merge sort algorithm we can sort all the $2^m$ values in $\mathcal{O}(2^m)$ steps.
Once groups are already $r$-sorted, discarding the dominated group minimals takes $\mathcal{O}(2^m)$ steps and Graham scan will now compute $\V$ in only
$\mathcal{O}(2^m)$ steps.  
Each look-up for a specific $\epsilon$ will now cost $\mathcal{O}(m).$   
 \end{remark} 

\subsection{Robustness to modeling errors}
\label{ss:epistemic}

Up till now, we have only considered aleatoric uncertainty,
taking on faith the correctness of input data: the probability distribution over the set of possible targets and their exact locations
were always assumed to be known. 
We next discuss the robustness with respect to epistemic (or ``systemic'') uncertainty caused by modeling errors.

\subsubsection{Robustness to errors in $\hat{\mathbf{\textit{p}}}$}
\label{sss:dro}
In practice, the probability distribution $\bphat$ over the set of targets will be usually defined using historical data.  
So, the true distribution would typically be different from such nominal/assumed $\bphat$.  
\textit{Distributionally robust} (DR) optimization deals with this uncertainty by defining an \textit{ambiguity set} of possible distributions, 
and letting the opponent/Nature select the worst among them while we are minimizing the expected cost \cite{delage2010distributionally}.
We will use the standard {\em total variation distance}, which for probability distributions $\bphat = (\phat_1,\cdots,\phat_m)$ and $\bptilde = (\ptilde_1,\cdots,\ptilde_m)$
over a finite set $\Xhat$ can be conveniently computed as
\begin{equation}
W(\bphat,\bptilde) = 1 - \sum_{i=1}^m \min(\phat_i, \ptilde_i) = \frac{1}{2} \|\bphat - \bptilde\|_1,
\label{eq:total_var_distance}
\end{equation}
coinciding with a Wasserstein-1 distance 
if we use a discrete metric 
$d(\xhat_i, \xhat_j) = 1 - \delta_{ij}$ on the set of targets \cite{villani2008optimal}.

We will refer to our assumed pdf $\bphat$ as the \textit{nominal distribution} and will define 
the ambiguity set to be a closed ball  
$B_\delta(\bphat) = \{ \bptilde: W(\bphat,\bptilde) \leq \delta \}$.
If $\qtilde(\x,\bptilde)$ denotes the expected time-to-target corresponding to $\bptilde$
and $\qtilde_\delta(\x)$ is the worst $\qtilde(\x,\bptilde)$ 
among all $\bptilde$'s
in $B_\delta(\bphat),$
the DR approach prescribes traveling to a waypoint $\bstilde_\delta$ that minimizes
that worst outcome.  In other words,
\[
\qtilde(\x,\bptilde) = \sum_{i=1}^m \ptilde_i u_i(\x),
\qquad
\qtilde_\delta(\x) = \max_{\bptilde \in B_\delta(\bphat)} \qtilde(\x,\bptilde),
\qquad
\bstilde_\delta \in \argmin_{\x \in \omegat} \, \qtilde_\delta(\x).
\]

Our explicit representation of $W(\bphat,\bptilde)$ makes computing $\qtilde_\delta(\x)$ particularly easy. 
Suppose that $u_i(\x)$'s are sorted in ascending order; i.e., $u_1(\x) \leq u_2(\x) \leq \cdots \leq u_m(\x)$.
The maximizer $\bptilde_\delta = (\ptilde_{\delta,1},\cdots,\ptilde_{\delta,m})$  will increase the probability of the ``hardest target'' $\xhat_m$  
while compensating with a corresponding decrease in probabilities of the ``easiest targets'' ($\xhat_1, \xhat_2, ...$). 
The first half of this plan can be implemented by setting $\ptilde_{\delta,m} = \min\{1, \phat_m + \delta\}.$
If $\ptilde_{\delta,m} = 1$ this immediately implies $\ptilde_{\delta,i} = 0$ for $i<m.$
Otherwise, this increase is compensated by starting with $\ptilde_{\delta,1} = \max\{0, \phat_1 - \delta\}$
and continuing (e.g., with $\ptilde_{\delta,2} = \max\{0, \phat_2 - \delta + (\phat_1 - \ptilde_{\delta,1})\}$)
until we reach a total reduction of $\delta$.
Since the difference between $u_i(\x)$'s is at most $[u_m(\x) - u_1(\x)]$, the DR loss can simply be bounded as $\qtilde_\delta(\x) \leq q(\x) + \delta [u_m(\x) - u_1(\x)]$.
This reveals that, for small $\delta$, minimizing $\qtilde_\delta$ is approximately equivalent to minimizing $q(\x)$.
Obtaining $\bstilde_\delta$ mainly involves computing $u_1, \cdots, u_m$, which takes $O(mM\log M)$ operations, and sorting $\phat_1,\cdots,\phat_m$ , which takes $O(m\log m)$ operations.

\begin{figure}[h]
\centering
$
\begin{array}{ccc}
\includegraphics[width = 0.3\textwidth]{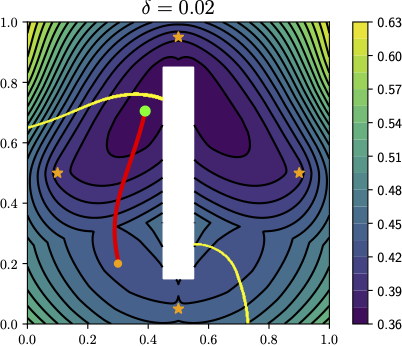} &
\includegraphics[width = 0.3\textwidth]{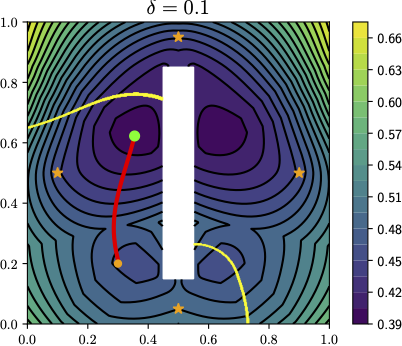} &
\includegraphics[width = 0.32\textwidth]{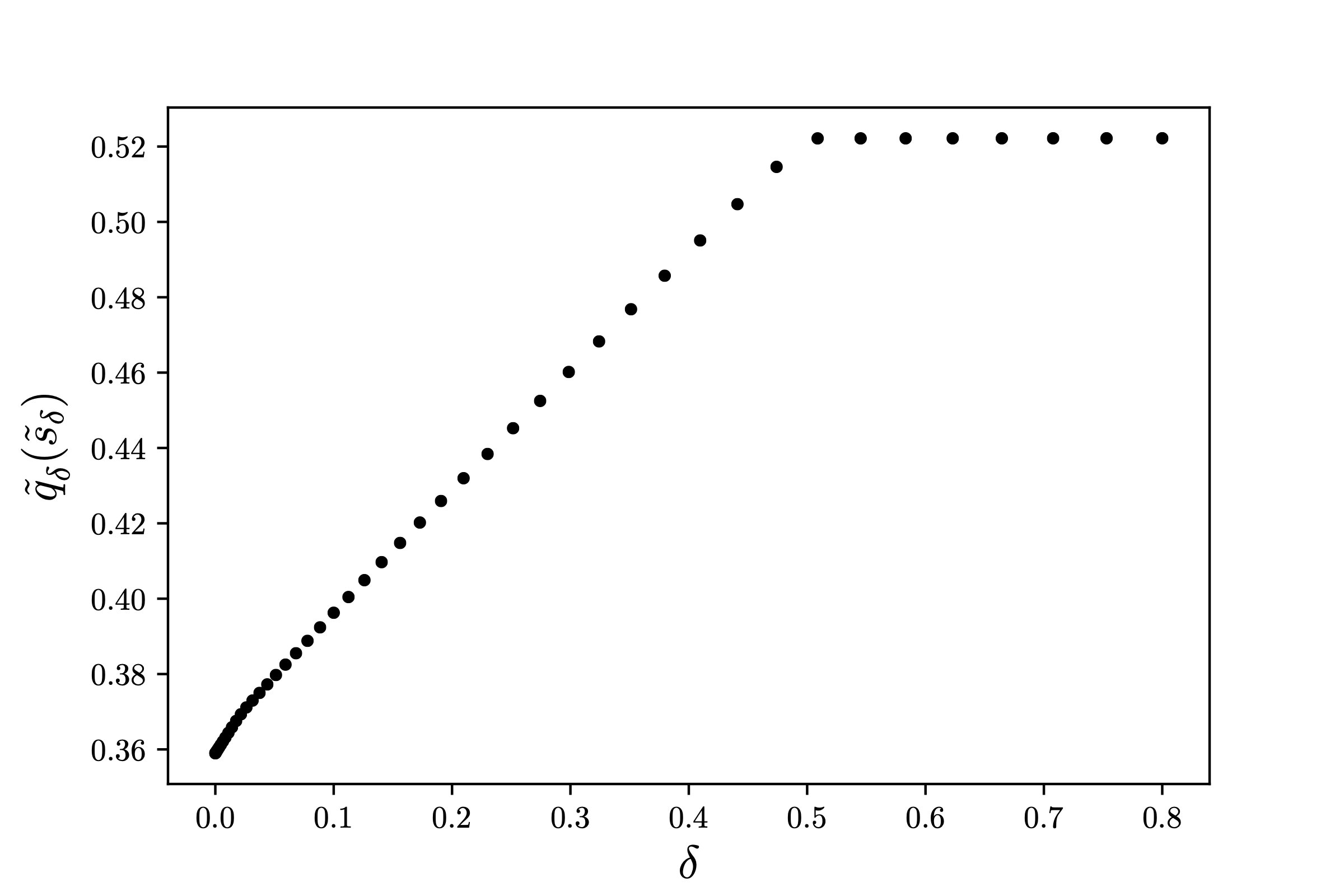}\\
(A) & (B) & (C)
\end{array}
$
\caption{DR optimization. 
(A) \& (B): $\delta = 0.02, 0.1$ respectively. The red curve is time-optimal to a waypoint $\bstilde_\delta$(lime dot). Level sets of $\qtilde_\delta$ are shown in the background.
(C): $\qtilde_\delta(\bstilde_\delta)$ computed using 40 different $\delta$ values between $0$ and $0.8$.
}
\label{fig:distribution_robust}
\end{figure}

In Figure \ref{fig:distribution_robust} we show the results of the DR approach for the fixed $T=0.4$ example of section \ref{sec:FixedT}.
For $\delta = 0$, there is no ambiguity in distribution and we are left with $\qtilde_\delta=q,$ the case presented in Figure \ref{fig:FixedT}.
For $\delta \geq (1 - \min_i \phat_i),$ we know that $\bptilde_\delta$ will assign probability 1 to the worst target, and $\qtilde_\delta = \qbar.$
(Note that this can happen even for smaller $\delta$'s; e.g., when 
the probability of the most difficult target is higher.
See Figure \ref{fig:distribution_robust}C.)
But for the intermediate $\delta$ values, $\qtilde_\delta$ continues to change and the corresponding waypoint $\bstilde_\delta$ will be different from $\sa$ and $\sw$.

It is also worth asking how good is the DR approach for balancing the worst-average case performance with respect to our nominal distribution $\bphat$.
The green dots in Figure \ref{f:worst_expected}B correspond to the same waypoints used in Figure \ref{fig:distribution_robust}.  
They clearly show that the DR waypoints are better in this sense than the risk-sensitive waypoints (defined in section \ref{ss:risk_sensitive}), but 
certainly are not Pareto optimal.  
An additional example in Appendix \ref{sec:bad_dr_example} shows that the DR approximation can also miss most of the Pareto frontier.

\subsubsection{Coarsening the target set}
All problems considered so far assumed that the target set is finite.
We now suppose that the target $\xhat$ is a random variable taking values in $\domain$, with a general probability measure $\measure$.
The expected minimum time-to-target is
\begin{equation}
\xi(\x) = \E_{\xhat} u(\x;\xhat) = \int\limits_{\domain} u(\x;\y)\, d\measure(\y),
\end{equation}
where $u(\x;\y)$ is obtained by solving \eqref{eq:Eik} with $\Gamma = \{\y\}$. 
If our goal is to achieve the average-case optimality with a fixed target identification time $T$, we should select a waypoint
\[
\sa_{\measure} \in \argmin\limits_{
\x \in \domainT(\source) } \xi(\x).
\]

However, for the sake of computational efficiency, we could instead use a domain decomposition $\domain = \bigcup\limits_{i = 1}^m \Upsilon_i,$
with $m$ non-overlapping subdomains $\Upsilon_i$ and a single representative target $\xhat_i \in \Upsilon_i$ selected from each of them.
The approach from section \ref{sec:FixedT} can be now viewed as finding an approximately optimal waypoint for a {\em coarsened target set} $\Xhat$
with
\begin{equation}
q(\x) = \sum\limits_{i=1}^m \int\limits_{\Upsilon_i} u(\x;\xhat_i)\, d\measure(\y) = \sum\limits_{i=1}^m\ u(\x;\xhat_i) \phat_i, \qquad \phat_i = \int\limits_{\Upsilon_i} \, d\measure(\y).
\end{equation}
Once we reach the waypoint $\sa \in \argmin_{\domainT} q(\x)$ and the true target $\y \in \domain$ is revealed, we will proceed there directly.  This strategy yields
the average time to target $\xi(\sa),$ and the natural question is how much worse it is than the optimal $\xi(\sa_{\measure}).$

Suppose our speed of motion is bounded below by $f_l > 0$ and all subdomains are small enough so that  $\forall \y \in \Upsilon_i$
there is a path of length at most $h$ connecting $\y$ to $\xhat_i$ without leaving $\Upsilon_i.$  Then the triangle inequality yields
\begin{eqnarray}
\nonumber
\left| \xi(\x) - q(\x) \right| 
& \leq & \sum\limits_{i=1}^m \int\limits_{\Upsilon_i} \left|u(\x;\y) - u(\x;\xhat_i)\right|\, d\measure(\y)
\; \leq \; \sum\limits_{i=1}^m \int\limits_{\Upsilon_i} u(\y; \xhat_i) \, d\measure(\y)\\
& \leq & \sum\limits_{i=1}^m \int\limits_{\Upsilon_i} (h/f_l) \, d\measure(\y) 
\; = \; h/f_l.
\label{eq:continuous_approx}
\end{eqnarray}

\noindent
Recalling the $q$-optimality of $\sa$ and using \eqref{eq:continuous_approx} twice (with $\x = \sa_{\measure}$ and then with $\x = \sa$),
we see that
\begin{equation}
\label{eq:continuous_estimate}
\xi(\sa_{\measure}) \geq q(\sa_{\measure}) - h/f_l \geq q(\sa) - h/f_l \geq \xi(\sa) - 2h/f_l,
\end{equation}
i.e., the suboptimality of  $\sa$ is bound by $2h/f_l$.


\section{Conclusions}
\label{sec:Conclusions}

\vspace{.5cm}
Our focus has been on the challenge of optimal and robust path-planning under initial uncertainty.
We have considered three different models for the certainty time $T$ 
and derived
efficient numerical methods for optimizing the average case performance in each case.  For the fixed $T$ and exponentially distributed $T$ we have also introduced a number of robust path-planning techniques, encoding different approaches to balancing the average case performance against the particularly bad/unlucky scenarios.  Under the assumption that the running cost $K$ and the certainty time $T$ are constant, the problem reduces to a choice among $T$-reachable waypoints, which allowed us to evaluate several robust planning techniques that are often prohibitively expensive in the general setting.  E.g., we have introduced a polynomial time algorithm for optimizing the expectation under probabilistic constraints on bad outcomes.  Risk-sensitive and distributionally-robust optimization are often viewed as closely approximating the optimal tradeoffs between the average and worst\change{-}case performance. But we have shown that in our setting, the resulting waypoints can be far from Pareto optimal.  Extending these robust planning algorithms to the case of general random $T$ is a challenge that we hope to address in the future.

While most of our presentation is focused on target uncertainty, a similar approach can be employed with initial uncertainty on cost or dynamics (on yet to be explored part of the domain);
see Appendix \ref{sec:storm_fronts} for a representative example.  For the sake of simplicity, our examples were always isotropic and relied on Fast Marching Method with first-order accurate upwind finite difference schemes.  But similar techniques would be easy to implement with Fast Marching-style methods developed for higher-order accurate discretizations (e.g., \cite{sethian2000fast, sethian1999fast, ahmed2011third})
and more general anisotropic problems (e.g., \cite{sethian2001ordered, alton2012ordered, dahiya2018ordered, mirebeau2014efficient, mirebeau2019riemannian}).

Beyond path planning, we believe that our approach is also promising for treating structured uncertainty in general control problems 
and Markov Decision Processes \cite{FlemingSoner}.  
For instance, techniques similar to those described in our section \ref{sec:DiscreteRandomT} have recently been applied to 
determining optimal car braking policy under traffic signal duration uncertainty
 \cite{gaspard2022optimal}.
 
Another useful future direction is to adapt our techniques to the setting of Mean Field Games (MFG), where a large number of interacting agents plan their actions selfishly,
and the initial uncertainty can be viewed as a particularly simple model of (instantaneous) ``common noise'' affecting all agents.  E.g., evacuating a building under uncertainty on which
exit will be open at a later/known time $T$ is an example of this type considered in \cite{achdou2019mean, carmona2022convergence}.
If the agents are cooperative and centrally controlled, this can be handled in the framework or Mean Field Control \cite{lauriere2021numerics_review} or, for a smaller number of agents,
by solving a high-dimensional HJB equation (e.g., using neural networks \cite{nakamura2021adaptive, meng2022sympocnet} or tensor decomposition techniques \cite{dolgov2021tensor}).
As an alternative to the mean field approximations,
sequential path-planning  \cite{chen2015safe, robinson2018efficient} can be also used to alleviate the curse of dimensionality if 
the agents are cooperating and can be assigned priority.   
We believe that our models of initial uncertainty would be beneficial in each of these settings.

Throughout this paper, we have assumed that the uncertainty disappears instantaneously. It will be useful and challenging to extend our methods to general problems with delayed information acquisition.   
Example in Appendix \ref{sec:drones_example} is actually a small step in this direction, with target-uncertainty reduced gradually at several distinct times $T_i.$
A different generalization of this type will be 
to consider problems where the transition to certainty is still instantaneous, but our choice of pre-certainty control affects the probability distribution of $T$.
This situation arises naturally in many ``stochastic shortest path with recourse'' problems \cite{polychronopoulos1996stochastic}.

\vspace*{3mm}
\noindent
{\bf Acknowledgements:}
The authors would like to thank Antonio Farah, Qianli Song, and Tristan Reynoso for their help in initial testing of some of these methods in Summer REU program at Cornell University.
The authors are also grateful to Dmitriy Drusvyatskiy for sharing his insight on 
the number of waypoints needed in non-deterministic optimal policies when optimizing with probabilistic constraints.

\vspace*{3mm}
\noindent
{\bf Author Contributions}  All three authors contributed to the development of the paper and the underlying research.

\vspace*{3mm}
\noindent
{\bf Funding} This work was supported in part by the National
Science Foundation (grants DMS-1738010 and DMS-1645643),
Air Force Office of Scientific Research (FA9550-22-1-0528), and the last author's Simons Foundation Fellowship.

\vspace*{2mm}
\noindent
{\bf Data availability} No data associated in the manuscript.

\vspace*{4mm}
\noindent
{\large \bf Declarations:}

\vspace*{3mm}
\noindent
{\bf Conflict of interest}  The authors declare no competing interests.

\vspace*{3mm}
\noindent
{\bf Ethical Approval} Not applicable.

\appendix
\appendixpage


\section{Uncertainty in running cost: the storm fronts example}
\label{sec:storm_fronts}

Returning to the setting of a fixed certainty time $T$ in section \ref{sec:FixedT},
we now show how our approach can similarly handle
initial uncertainty arising due to anticipated 
changes in a global environment even if the target is perfectly known in advance.

Consider a challenge of risk-hedged air traffic routing under uncertain weather conditions.
Given a forecast of probable storm front trajectories, one option is to devise a good flight path and commit to it in advance all the way to the target.  
This is the key idea of the approach developed in  \cite{sadovsky2016risk}, but one can do much better by taking into account
the future re-routing once a more precise forecast becomes available.  What follows is an illustration of how this can be done in our framework under a 
simplifying assumption that the true position of the storm front will be revealed at a known time $T$.

\begin{figure}[h]
\centering
$
\begin{array}{ccc}
\includegraphics[width = 0.3\textwidth]{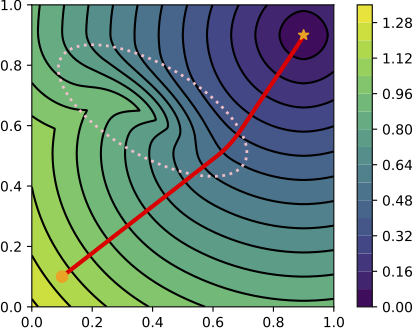} &
\includegraphics[width = 0.3\textwidth]{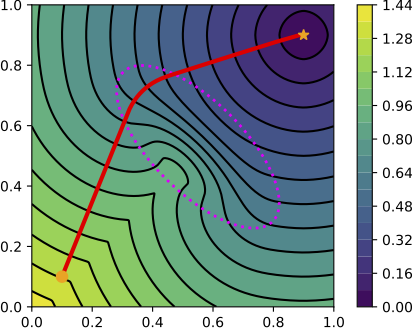} &
\includegraphics[width = 0.3\textwidth]{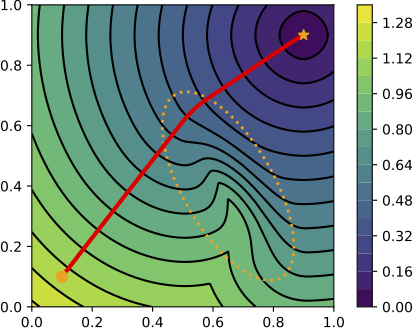}\\
(A) & (B) & (C)
\end{array}
$
\caption{A deterministic case of different storms with elliptic shapes. $f = 1, K = 1$ outside storm region, while $f = 1, K > 1$ inside. 
The source $\source$ is at $(0.1,0.1)$ while the target $\xhat$ is at $(0.9, 0.9)$. 
Each storm location results in a different running cost $K_i(\x)$ and value function $u_i(\x)$.
The red lines are cost-minimal trajectories and in the background are level sets of $u_1, u_2$ and $u_3$. 
}
\label{fig:storms}
\end{figure}

First, consider the deterministic case with an a priori known location of an elliptic storm region.  Our running cost $K$ is assumed to be identically one on most of the domain, but is higher inside the storm region.  In particular,  if the storm is centered at $\xtilde$ and the canonical equation for the storm region boundary is $(\x-\xtilde)^T A (\x-\xtilde) = 1$, we will assume that inside that region $K(\x) = 1 + \alpha \big(1 - (\x-\xtilde)^T A (\x-\xtilde) \big)^\gamma$ with $\alpha = 2$ and $\gamma = 2.5.$ (Note that this yields $K=1$ on the storm boundary, ensuring the continuity.)  For any fixed/known $\xtilde$ and $A$, finding an optimal path to the target $\xhat$ is simply a matter of solving an Eikonal equation and using gradient descent on the value function.   In Figure \ref{fig:storms} we illustrate this for three different storm locations (encoded by $(\xtilde, A)$). 
Unsurprisingly, optimal trajectories deviate from the quickest (straight line) path to the target to decrease the amount of time spent inside the storm.

Suppose that initially we only have a probability distribution $\bphat = (\phat_1,\phat_2,\phat_3)$ over these three storm regions (enumerated in the same order as in Figure \ref{fig:storms}) and the true position of the storm will be revealed at a later time $T$.  We will further assume that $T$ is small enough so that all three possible storm regions have no overlap with $\omegat;$ i.e., the plane can never reach any storm region before the ``storm-revealing time''. 
The optimal waypoint is found in the same way as in section \ref{sec:FixedT}.
I.e., the plane should follow the time-optimal path to 
\[
\sa \in \argmin\limits_{\x \in \domainT(\source)}  q(\x),
\qquad \qquad
\text{where } \;
 q(\x) \, = \, \E^{}_{\xtilde, A} \, \big[ u(\x) \big] \, = \, \sum_{i=1,2,3} \phat_i u_i(\x).
\]
In general, the time complexity of obtaining $\sa$ is $O(m_s M \log M)$ with $m_s$ being the total number of potential storm regions and $M$ the number of gridpoints.

\begin{figure}[!htb]
\centering
$
\begin{array}{cc}
\includegraphics[width = 0.45\textwidth]{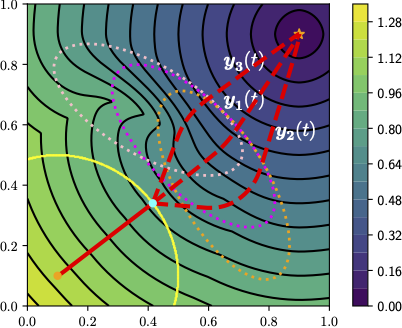} &
\includegraphics[width = 0.45\textwidth]{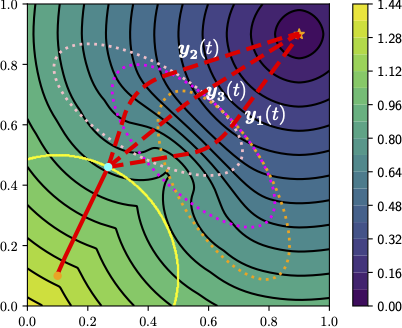}\\
(A) & (B)
\end{array}
$
\caption{Uncertain elliptic storms with different $\bphat$'s. Level sets of $q(\x)$ are shown in \change{the} background. The cyan dot is $\sa$, attaining minimal $q$ in $\omegat$, with $T = 0.4$. 
The dashed red lines labeled $\y_i(t)\, (i=1,2,3)$ are the optimal trajectories from $\sa$ to the target 
-- to be used once the actual storm location is discovered.
(A): $\bphat = (0.8, 0.1, 0.1)$. (B): $\bphat = (0.1, 0.8, 0.1)$. 
}
\label{fig:storm_expected}
\end{figure}

As we can see from Figure \ref{fig:storm_expected}, the optimal waypoint is heavily dependent on $\bphat$.
Of course, the above setting is highly simplified since in reality the weather forecasts would be updated numerous times or even continuously.
But the amount of uncertainty would normally still decrease with time and the main ideas of sections \ref{sec:FixedT}-\ref{sec:ExponentialRandomT} would be still applicable.


\section{Temporal changes in target probabilities}
\label{sec:drones_example}

The following is \change{an}  ``emergency rescue'' example where $p_j$'s are not independent of $\phat_i$'s.
The rescue vehicle starts at $\source$ and moves with isotropic speed $f$ in the domain with a rectangular obstacle, exactly as shown in Figure \ref{fig:determ}.  
The goal is to minimize the expected time until we rescue a subject known to be located at one of the 4 possible sites $\xhat_1,\cdots, \xhat_4$ with the a priori likelihood reflected by $\phat_1,\cdots, \phat_4.$ 
But unlike the previous examples, here we assume the availability of 3 lightweight aerial drones also launched from $\source$ and flying to the closest 3 sites ($\xhat_3, \xhat_4,$ and $\xhat_2$) along the straight lines (unconstrained by the ground rectangular obstacle).  For simplicity, we will assume that drones move with  speed $f_d=5/3$ and 
notify the main vehicle whether the rescue subject is
found as soon as they reach their respective destinations.
For the locations specified in Figure \ref{fig:determ},
these drones' arrival times will be 
$
T_1 = |\xhat_3 - \source| / f_d 
\, < \,
T_2 = |\xhat_4 - \source| / f_d 
\, < \,
T_3 = |\xhat_2 - \source| / f_d. 
$
If the first drone discovers that subject at $\xhat_3,$ the rescue vehicle follows the quickest path to it.  Otherwise, the probabilities of other sites is adjusted accordingly:  
$\phat^{new}_j = \phat_j / (1 -\phat_3)$ for all $j \neq 3$
and $\phat^{new}_3 = 0$.
The probabilities are similarly adjusted  at $T_2$ if $\xhat \neq \xhat_4.$ If the subject has not been located earlier, at the time $T_3$ we will discover whether he is at the last drone-visited site $\xhat_2$ or at the remaining site $\xhat_1.$  The question is how to optimally plan the rescue vehicle's path until the discovery time $T \in \{T_1, T_2, T_3\}.$ This can be handled by solving the same sequence of time-dependent PDEs \eqref{eq:tEik}.
The value functions $v_i(\x,t)$ will satisfy \eqref{eq:tEik} on $[T_{i-1},T_i], i = 1,2,3$ with $T_0 = 0$ and the terminal conditions  
\begin{align*}
v_3 \left(\x,T_3 \right) \; & = \; 
\frac{\phat_2}{\phat_2 + \phat_1} u_2 (\x) \, + \, 
\frac{\phat_1}{\phat_2 + \phat_1} u_1 (\x), \\
v_2(\x,T_2) \;  & = \; 
\frac{\phat_4}{\phat_1 + \phat_2 + \phat_4} u_4(\x) \, + \,\frac{\phat_1 + \phat_2}{\phat_1 + \phat_2 + \phat_4} v_3 \left(\x,T_2 \right), \\
v_1(\x,T_1) \; & = \; 
\phat_3 u_3(\x) \, + \, \left(\phat_1 + \phat_2 + \phat_4 \right) v_2 \left(\x,T_1 \right). \\
\end{align*}

\vspace*{-3mm}
\begin{figure}[h]
\label{fig:drones}
\centering
$
\arraycolsep=15pt\def\arraystretch{2}
\begin{array}{cc}
\includegraphics[width = 0.4\textwidth]{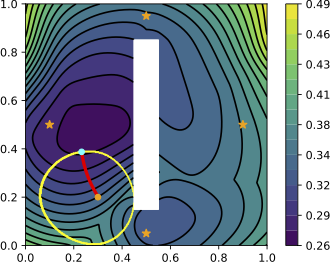} &
\includegraphics[width = 0.4\textwidth]{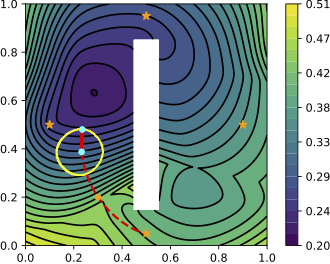}\\
(A) & (B)\\[3em]
\includegraphics[width = 0.4\textwidth]{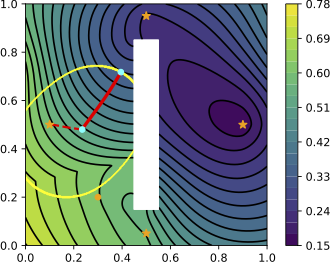} &
\includegraphics[width = 0.4\textwidth]{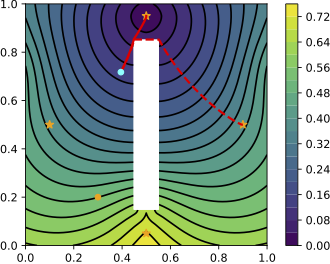}\\
(C) & (D) \\
\end{array}
$
\caption{The drone-assisted emergency rescue example with the same initial $\phat_i$'s as in Figure \ref{fig:FixedT}.
$(A)$: The optimal waypoint $\sa_1$ for  \change{the} time $T_1$ and the $\source \rightarrow \sa_1$ trajectory, superimposed on the level sets of $v_1(\x, T_1)$. The yellow curve shows the boundary of the  reachable region $\domain_{T_1}(\source)$.
$(B)$: If the first drone discovers the subject,
the vehicle follows a dashed red $\sa_1 \rightarrow \xhat_3$ trajectory.
Otherwise, it follows a solid red trajectory to the $T_2$-optimal waypoint $\sa_2$, shown on top of the level sets of $v_2(\x, T_2)$. The yellow curve is the boundary of $\domain_{T_2-T_1}(\sa_1)$.
$(C)$: If the subject is not discovered by the first drone but the second,
the vehicle follows a dashed red $\sa_2 \rightarrow \xhat_4$ trajectory.
Otherwise, it follows a solid red trajectory to the $T_3$-optimal waypoint $\sa_3$, shown on top of the level sets of $v_3(\x, T_3)$. The yellow curve is the boundary of $\domain_{T_3-T_2}(\sa_2)$.
$(D)$: By the time $T_3$, the drone either finds the subject at $\xhat_2$ and the vehicle follows the dashed red trajectory, or we conclude that the subject is at $\xhat_1$ and the vehicle follows the solid red path shown over the level sets of $u_1(\x)$. 
}
\end{figure}

\vspace*{3mm}

\FloatBarrier


\section{DR optimization missing a large part of $(q,\qbar)$ Pareto front}
\label{sec:bad_dr_example}

Figure \ref{f:worst_expected}B showed that risk-sensitive waypoints can remarkably deviate from $q-\qbar$ Pareto front. 
The distributionally-robust waypoints are 
much closer to Pareto optimality in that same figure.
However, the following example shows that often there is no DR-optimal $\bstilde_\delta$ satisfying $\qbar(\bstilde_\delta) = \qbar(\saw)$ for many $\saw$ values.
Thus, the DR results cannot approximate the corresponding parts of \change{the} Pareto front, which in this case is highly non-convex and even discontinuous. 
The DR optimal waypoints cluster into three groups, with large gaps between clusters.

\begin{figure}[h]
\centering
$
\begin{array}{ccc}
\includegraphics[width = 0.3\textwidth]{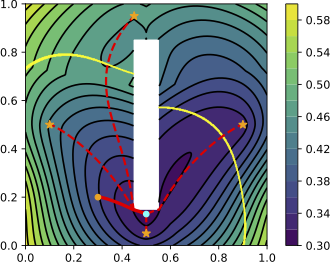} &
\includegraphics[width = 0.3\textwidth]{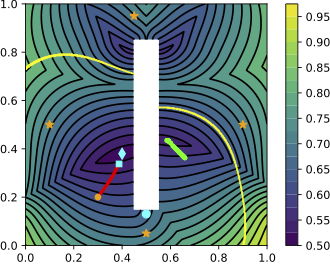} &
\includegraphics[width = 0.32\textwidth]{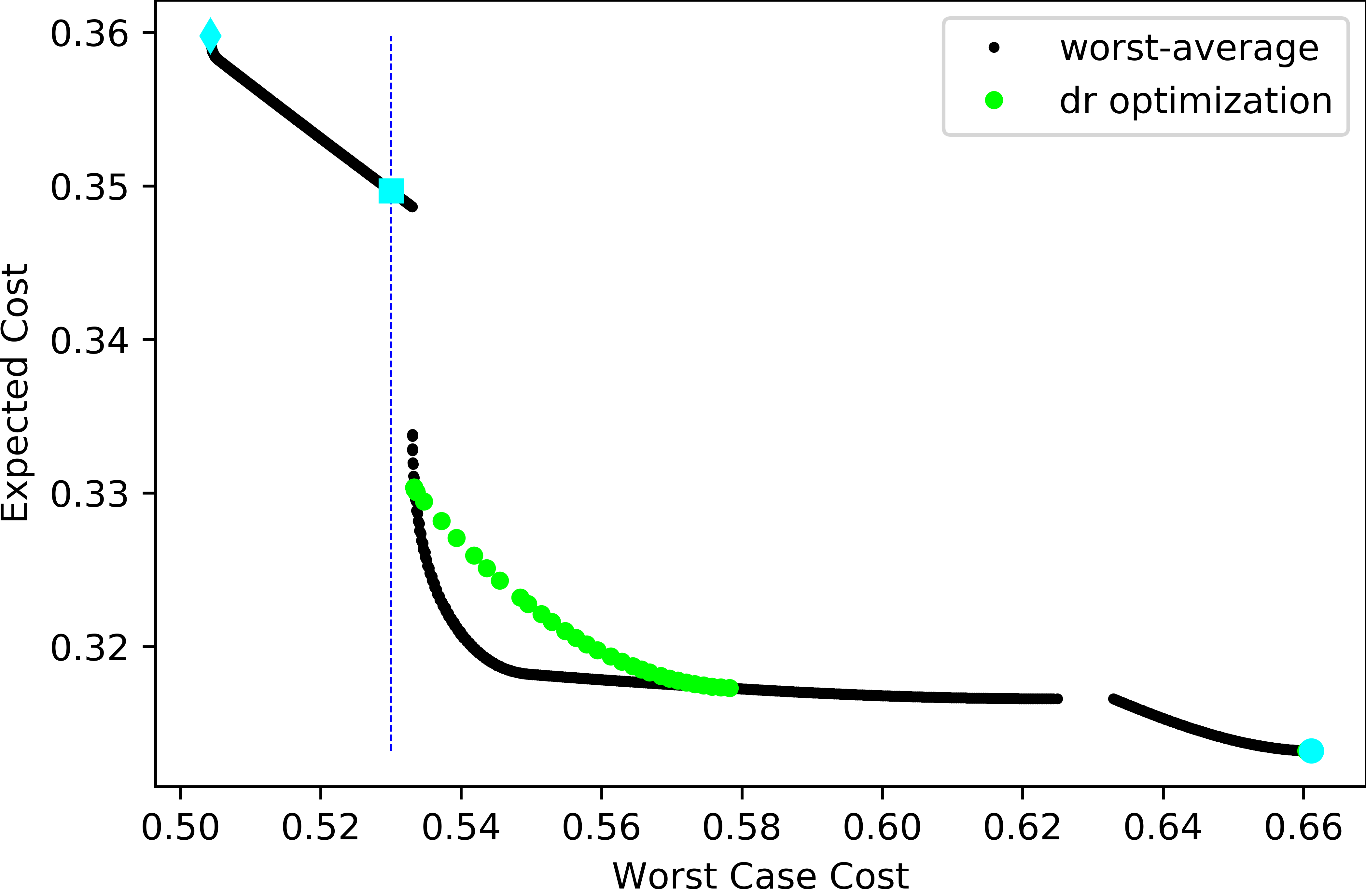}\\
(A) & (B) & (C)
\end{array}
$
\caption{
(A): Expectation-optimal waypoint and post-certainty paths to each target. 
The corresponding target probabilities (clockwise, starting at the top) are $\phat = (0.17, 0.35, 0.3, 0.18)$.
The symmetry of target locations is broken as the top target ($\xhat_1$) is moved slightly left to $(0.45, 0.95)$. 
The speed function is also slightly different: $f= 1.4 - 0.6cos(2\pi x)sin(2\pi y)$.
Level sets of $q(\x)$ are shown in the background.
(B): Three different optimal waypoints ($\sa,  \sw, \saw$) (cyan dot, diamond and square) corresponding to average, worst-case optimality and minimizing average with a hard constraint $C = 0.53$.
A hundred distributionally-robust waypoints $\bstilde_\delta$'s are also shown as green dots with equally spaced $\delta$s between 0 and 1.
Level sets of $\qbar(\x)$ are shown in the background.
(C): $q-\qbar$ Pareto front and all the waypoints from (B). 
}
\label{fig:bad_dr_example}
\end{figure}

\FloatBarrier

\bibliographystyle{siam}
\bibliography{Ref}

\end{document}